\documentclass{siamltex}

\usepackage{amsmath,amsfonts,graphicx,color,multirow}
\usepackage{subfig,booktabs,array,hyphenat,amsbsy}
\usepackage[breaklinks=true,pdfborder={0 0 0 0}]{hyperref}

\newcolumntype{C}[1]{>{\centering\let\newline\\\arraybackslash\hspace{0pt}}m{#1}}

\newcommand{\porofluid}{poroelastic\hyp fluid}

\newcommand{\clawpack}{{\sc clawpack}}

\newcommand{\refcartpaper}{\cite{gil-ou-rjl:poro-2d-cartesian}}
\newcommand{\refmappaper}{\cite{gil-ou:poro-2d-mapped}}

\newcommand{\beq}{\begin{eqnarray}}
\newcommand{\eeq}{\end{eqnarray}}
\newcommand{\beqs}{\begin{eqnarray*}}
\newcommand{\eeqs}{\end{eqnarray*}}
\newcommand{\bC}{\mathbf{C}}

\newcommand{\btau}{\mbox{\boldmath $\tau$}}

\newcommand{\bA}{\mbox{$\mathbf{A}$}}
\newcommand{\bB}{\mbox{$\mathbf{B}$}}
\newcommand{\bD}{\mbox{$\mathbf{D}$}}
\newcommand{\bE}{\mbox{$\mathbf{E}$}}
\newcommand{\bG}{\mbox{$\mathbf{G}$}}

\newcommand{\balpha}{\mbox{\boldmath $\alpha$}}

\newcommand{\bQ}{\mbox{$\mathbf{Q}$}}

\newcommand{\B}{\mathbf}

\newcommand{\vect}[1]{\mathbf{#1}}
\DeclareMathOperator{\divg}{div}

\DeclareMathOperator{\real}{Re}
\DeclareMathOperator{\imag}{Im}
\newcommand*\Bell{\ensuremath{\boldsymbol\ell}}

\author{Grady\ I.\ Lemoine\footnotemark[2]}

\begin{document}

\title{Three-Dimensional Mapped-Grid Finite Volume Modeling of Poroelastic-Fluid Wave Propagation\footnotemark[1]}
\maketitle

\renewcommand{\thefootnote}{\fnsymbol{footnote}}

\footnotetext[1]{Much of this work comes from Section 2.2 and Chapters
  8 and 9 of~\cite{lemoine:thesis}.  This is the first
  widely-distributed version.}

\footnotetext[2]{Department of Applied
  Mathematics, University of Washington, Seattle, WA 98195
  (gl@uw.edu).  Supported in part by NIH grant
  5R01AR53652-2 and NSF grants DMS-0914942 and DMS-1216732.}

\begin{abstract}
This paper extends the author's previous two-dimensional work with Ou
and LeVeque to high-resolution finite volume modeling of systems of
fluids and poroelastic media in three dimensions, using logically
rectangular mapped grids.  A
method is described for calculating consistent cell face areas and
normal vectors for a finite volume method on a general non-rectilinear
hexahedral grid.  A novel limiting algorithm is also developed to cope
with difficulties encountered in implementing high-resolution finite
volume methods for anisotropic media on non-rectilinear grids; the new
limiting approach is compatible with any limiter function, and
typically reduces solution error even in situations where it is not
necessary for correct functioning of the numerical method.
Dimensional splitting is used to reduce the computational cost of the
solution.  The code implementing the three-dimensional algorithms is
verified against known plane wave solutions, with particular attention
to the performance of the new limiter algorithm in comparison to the
classical one.  An acoustic wave in brine striking an uneven bed of
orthotropic layered sandstone is also simulated in order to demonstrate the capabilities
of the simulation code.
\end{abstract}



\section{Introduction}
\label{sec:intro}

Biot poroelasticity theory is a homogenization technique for modeling
the mechanics of a fluid-saturated porous solid.  It was developed in
the period from the 1930s to the 1960s for problems in soil and rock
mechanics~\cite{biot:56-1, biot:56-2, biot:62}, but has also found use
in {\em in vivo} bone~\cite{cowin:bone-poro,
  cowin-cardoso:fabric-2010, gilbert-guyenne-ou:bone-2012} and
underwater acoustics~\cite{bgwx:2004, gilbert-lin:1997,
  gilbert-ou:2003}.

In Biot theory, the solid part of the medium (termed the {\em matrix}
or {\em skeleton}) is modeled using linear elasticity, while the fluid
is treated using compressible linearized fluid dynamics; Darcy's law
is used to model the aggregate motion of the fluid through the matrix.
The interaction of the fluid and solid gives rise to three different
types of waves: fast P waves, which are similar to the P waves of
elastodynamics; shear waves similar to elastodynamic S waves; and slow
P waves, which produce behavior not found in simpler types of media.
The interaction of the fluid with the solid part of the medium is
important to the behavior of these waves --- for the fast P and shear
waves, there is little motion of the fluid with respect to the solid,
but the slow P waves show relatively large amounts of fluid motion.
The viscosity of the pore fluid thus causes light damping and
dispersion of the first two wave types, but strong slow P wave damping
and a substantial variation of slow P phase velocity with frequency.
Carcione provides an excellent treatment of poroelasticity theory in
Chapter 7 of his book~\cite{carcione:wave-book}.

A variety of methods have been used to model poroelasticity
numerically, including finite difference and
pseudospectral~\cite{chiavassa-lombard:acou-poro-ccp,
  dai-vafidis-kanasewich:poro-vel-stress, garg-nayfeh-good:poro-1974,
  mikhailenko:poro-1985}, finite
element~\cite{buchanan-gilbert-khashanah:bone-params-lofreq-2004,
  santos-orena:poro-fem}, boundary
element~\cite{attenborough-berry-chen:scattering}, spectral
element~\cite{degrande-deroeck:poro-spectral-freq,
  morency-tromp:poro-spectral-time}, discontinuous
Galerkin~\cite{delapuente-dumbser-kaser-igel:poro-dg}, and finite
volume~\cite{gil-ou-rjl:poro-2d-cartesian, naumovich:fvm-3d-staggered}
methods.  Semi-analytical methods have been used as well, both in
forward problems for their own sake~\cite{detournay-cheng:borehole}
and as the forward solution component of an inversion
scheme~\cite{buchanan-gilbert:bone-params-hifreq-1-2007,
  buchanan-gilbert:bone-params-hifreq-2-2007}.

Three-dimensional poroelasticity simulation has become somewhat common
in recent years due to the increasingly powerful computers available.
Recent works on three-dimensional computational poroelasticity on
regular grids include those of
Naumovich~\cite{naumovich:fvm-3d-staggered}, who used a staggered-grid
finite volume method on regular, rectilinear grids for isotropic
media, and Aldridge et al.~\cite{aldridge-symons-bartel:fd-3d}, who
used a staggered-grid finite difference approach.  Three-dimensional
work capable of using irregular grids includes that of de la Puente et
al.~\cite{delapuente-dumbser-kaser-igel:poro-dg}, who employed a
discontinuous Galerkin method on triangular and tetrahedral meshes,
and the spectral element work of Morency and
Tromp~\cite{morency-tromp:poro-spectral-time}.

This paper extends the previous work of Lemoine, Ou, and
LeVeque~\cite{gil-ou:poro-2d-mapped, gil-ou-rjl:poro-2d-cartesian} to systems of orthotropic
poroelastic and fluid media in three dimensions, modeled using
logically rectangular mapped grids.  Section \ref{sec:theory} develops
a first-order linear system of PDEs modeling low-frequency Biot
poroelasticity theory, and repeats the interface conditions derived
in~\refmappaper{} for convenient reference.  Following this, section
\ref{sec:mapped-fvm-3d} extends the numerical methods of the previous
papers to the three-dimensional system.  While most of the extension
process is straightforward, some problems occur in 3D that have no
counterpart in 2D; section \ref{sec:mapped-grids-3d} discusses a
complication that arises when implementing a finite volume method on a
non-rectilinear grid of hexahedral cells, while section
\ref{sec:elimiter} formulates a new wave strength ratio for wave
limiting in order to circumvent problems with consistent shear wave
identification for mapped grids on orthotropic media.  The
simulation code implementing these numerical methods is then verified
against known plane wave solutions in Section \ref{sec:results}, with
special attention paid to the behavior of the new limiter algorithm,
and results are shown for a poroelastic/fluid demonstration problem
showing the ability of the code to model complex systems on mapped
grids.

\section{Governing equations in three dimensions}
\label{sec:theory}

\subsection{Stress rate-velocity relations}

Equations (7.131) and (7.133) of Carcione~\cite{carcione:wave-book} give the
stress-strain relation for an anisotropic poroelastic material in an
orthogonal set of axes labeled 1, 2, and 3.  Using the summation
convention for repeated indices, these equations are
\begin{equation} \label{eq:aniso-stress-strain}
    p = M (\zeta - \alpha_I e_I), \quad
    \tau_I = c_{IJ}^u e_J - M \alpha_I \zeta.
\end{equation}
The quantities in this system are defined as follows:
\begin{itemize}
\item $p$ is the fluid pressure
\item $\zeta$ is the variation of fluid content, $\zeta = -\divg(\phi(
  \vect{u}_f - \vect{u}_m))$
\item $\phi$ is the porosity of the material
\item $\vect{u}_f$ and $\vect{u}_m$ are the displacements of the fluid
  and solid, respectively, from their stress-free configurations
\item $e_I$ is the $I$'th component of engineering strain,
  $e = \left(\begin{smallmatrix} \epsilon_{11} & \epsilon_{22} & \epsilon_{33} &
  \gamma_{23} & \gamma_{13} & \gamma_{12} \end{smallmatrix} \right)^T$.  Note that
  the engineering strains $\gamma_{ij}$ are twice the
  tensor strains $\epsilon_{ij}$, $i \ne j$.
\item $\tau_I$ is the $I$'th component of the total stress in the
  material, taken in the same order as the strain
\item $c_{IJ}^u = c_{IJ} + \alpha_I\alpha_J M$ is the undrained
  elastic stiffness tensor of the matrix
\item $c_{IJ}$ is the drained elastic stiffness tensor of the matrix
\item $\alpha_I$ is the $I$'th effective stress coefficient, given by
  $\alpha_I = 1 - \frac{1}{3K_s} \sum_{J=1}^3 c_{IJ}$
\item $K_s$ is the bulk modulus of the matrix material
\item $M$ is a parameter related to the bulk compressibility of the
  medium,
  \begin{equation}
    M = K_s \left( \left( 1 - \frac{K^*}{K_s} \right) - \phi \left( 1
    - \frac{K_s}{K_f} \right) \right)^{-1}
  \end{equation}
\item $K_f$ is the bulk modulus of the fluid
\item $K^*$ is another bulk stiffness coefficient, $K^* = \frac{1}{9}
  \sum_{I=1}^3 \sum_{J=1}^3 c_{IJ}$
\end{itemize}

To begin building a first-order velocity-stress
system, note the following relations between velocities and strain
rates (presuming space and time derivatives can be interchanged):
\begin{gather} \label{eq:strain-rate-velocity}
  \begin{aligned}
    \partial_t e_1 &= v_{1,1}, \quad
    &\partial_t e_2 &= v_{2,2}, \quad
    &\partial_t e_3 &= v_{3,3},\\
    \partial_t e_4 &= v_{2,3} + v_{3,2}, \quad
    &\partial_t e_5 &= v_{1,3} + v_{3,1}, \quad
    &\partial_t e_6 &= v_{1,2} + v_{2,1},\\
  \end{aligned}\\
  \partial_t \zeta = -q_{1,1} - q_{2,2} - q_{3,3}.
\end{gather}
Here $\vect{v}$ is the velocity of the matrix relative to an inertial
frame, and $\vect{q}$ is the flow rate of the fluid relative to the
matrix (the porosity $\phi$ times the aggregate velocity of the fluid
relative to the matrix).  Subscript indices before a comma represent
components, while those after
a comma represent differentiation.  Differentiating $\eqref{eq:aniso-stress-strain}$ with
respect to time, and defining the vectors of stresses and velocities as
$\bQ_s = \left(\begin{smallmatrix} \tau_{11} & \tau_{22} & \tau_{33} &
  \tau_{23} & \tau_{13} & \tau_{12} & p \end{smallmatrix} \right)^T$ and $\bQ_v
= \left(\begin{smallmatrix} v_1 & v_2 & v_3 & q_1 & q_2 &
  q_3 \end{smallmatrix} \right)^T$,
results in a system relating $\partial_t \bQ_s$ to the gradients of the
velocities:
\begin{equation} \label{eq:Asv-system}
  \partial_t \bQ_s + \bA_{sv} \bQ_{v,1} + \bB_{sv} \bQ_{v,2} +
  \bC_{sv} \bQ_{v,3} = 0.
\end{equation}

Rather than give the matrices $\bA_{sv}$, $\bB_{sv}$, and $\bC_{sv}$
individually, it is more convenient to write the matrix
$\breve{\bA}_{sv} = n_1 \bA_{sv} + n_2 \bB_{sv} + n_3 \bC_{sv}$ that
is the coefficient of the directional derivative of $\bQ_v$ in the
$(n_1, n_2, n_3)$ direction.  If the medium is orthotropic, and the
1-2-3 axes are its principal axes, then
\begin{equation} \label{eq:Asv}
  \breve{\bA}_{sv} = - \begin{pmatrix}
    n_1 c_{11}^u & n_2 c_{12}^u & n_3 c_{13}^u & n_1 \alpha_1 M & n_2 \alpha_1 M & n_3 \alpha_1 M \\
    n_1 c_{12}^u & n_2 c_{22}^u & n_3 c_{23}^u & n_1 \alpha_2 M & n_2 \alpha_2 M & n_3 \alpha_2 M \\
    n_1 c_{13}^u & n_2 c_{23}^u & n_3 c_{33}^u & n_1 \alpha_3 M & n_2 \alpha_3 M & n_3 \alpha_3 M \\
    0 & n_3 c_{44} & n_2 c_{44} & 0 & 0 & 0 \\
    n_3 c_{55} & 0 & n_1 c_{55} & 0 & 0 & 0 \\
    n_2 c_{66} & n_1 c_{66} & 0 & 0 & 0 & 0 \\
    -n_1 M \alpha_1 & -n_2 M \alpha_2 & -n_3 M \alpha_3 & -n_1 M & -n_2 M & -n_3 M
  \end{pmatrix}.
\end{equation}

\subsection{Equations of motion}

System \eqref{eq:Asv-system} does not yet provide a
closed set of equations that can be used to describe the dynamics of
the poroelastic medium.  Equations of motion are still required that
relate accelerations to gradients of stress.  Equations (7.255) and
(7.256) of~\cite{carcione:wave-book} provide the key.  If the medium
is orthotropic and the 1-2-3 axes are its principal axes, 
they relate accelerations to stress gradients:
\begin{equation} \label{eq:eqns-of-motion}
  \begin{aligned}
    \tau_{ij,j} &= \rho \partial_t^2 u_{mi} + \rho_f \partial_t^2 w_i \\
    -p_{,i} &= \rho_f \partial_t^2 u_{mi} + m_i \partial_t^2 w_i +
    \frac{\eta}{\kappa_i} \partial_t w_i.
  \end{aligned}
\end{equation}
Here $i$ ranges from 1 to 3; there is a sum over $j$ in the first
equation, but no sum over $i$ in the second.  The new variables in
these equations are as follows:
\begin{itemize}
\item $\rho$ is the bulk density of the medium, $\rho = (1 - \phi)
  \rho_s + \phi \rho_f$
\item $\rho_s$ is the density of the matrix material
\item $\rho_f$ is the density of the fluid
\item $\vect{w}$ is the displacement of the fluid relative to the
  matrix, scaled by the porosity.  The rate of change of $\vect{w}$ is
  $\partial_t\vect{w} = \vect{q}$.
\item $m_i$ is the fluid inertia along axis $i$, $m_i = \rho_f
  T_i/\phi$
\item $T_i$ is the tortuosity of the matrix along axis $i$, defined as
  the factor by which the kinetic energy of the fluid must be higher
  than its density would indicate for straight-line motion, in order
  to have a given bulk velocity along that axis
\end{itemize}

For each $i$, \eqref{eq:eqns-of-motion} is a system of two equations
in two unknowns.  Noting that $\partial_t u_{mi} = v_i$, this system becomes
\begin{equation} \label{eq:eom-intermed}
  \begin{aligned}
    \rho \partial_t v_i + \rho_f \partial_t q_i &= \tau_{ij,j} \\
    \rho_f \partial_t v_i + m_i \partial_t q_i &= -p_{,i} -
    \frac{\eta}{\kappa_i} q_i.
  \end{aligned}
\end{equation}
Solving with Cramer's Rule results in
\begin{equation} \label{eq:eom-accels}
  \begin{aligned}
    \partial_t v_i &= \frac{m_i}{\Delta_i} \tau_{ij,j} +
    \frac{\rho_f}{\Delta_i} p_{,i} + \frac{\rho_f \eta}{\Delta_i
      \kappa_i} q_i\\
    \partial_t q_i &= -\frac{\rho_f}{\Delta_i} \tau_{ij,j} -
    \frac{\rho}{\Delta_i} p_{,i} - \frac{\rho \eta}{\Delta_i \kappa_i} q_i,
  \end{aligned}
\end{equation}
where $\Delta_i := \rho m_i - \rho_f^2$.

It is now possible to write a linear system relating the rates of
change of the velocities to the gradients of stress, of the form
\begin{equation} \label{eq:Avs-system}
  \partial_t \bQ_v + \bA_{vs} \bQ_{s,1} + \bB_{vs} \bQ_{s,2} + \bC_{vs} \bQ_{s,3} =
  \bD_v \bQ_v.
\end{equation}
Again, it is more convenient to provide $\breve{\bA}_{vs} = n_1
\bA_{vs} + n_2 \bB_{vs} + n_3 \bC_{vs}$, rather than the individual
matrices of \eqref{eq:Avs-system}:
\begin{equation} \label{eq:Avs}
  \breve{\bA}_{vs} = - \begin{pmatrix}
    n_1 \frac{m_1}{\Delta_1} & 0 & 0 & 0 & n_3 \frac{m_1}{\Delta_1} & n_2 \frac{m_1}{\Delta_1} & n_1 \frac{\rho_f}{\Delta_1} \\
    0 & n_2 \frac{m_2}{\Delta_2} & 0 & n_3 \frac{m_2}{\Delta_2} & 0 & n_1 \frac{m_2}{\Delta_2} & n_2 \frac{\rho_f}{\Delta_2} \\
    0 & 0 & n_3 \frac{m_3}{\Delta_3} & n_2 \frac{m_3}{\Delta_3} & n_1 \frac{m_3}{\Delta_3} & 0 & n_3 \frac{\rho_f}{\Delta_3} \\
    -n_1 \frac{\rho_f}{\Delta_1} & 0 & 0 & 0 & -n_3 \frac{\rho_f}{\Delta_1} & -n_2 \frac{\rho_f}{\Delta_1} & -n_1 \frac{\rho}{\Delta_1} \\
    0 & -n_2 \frac{\rho_f}{\Delta_2} & 0 & -n_3 \frac{\rho_f}{\Delta_2} & 0 & -n_1 \frac{\rho_f}{\Delta_2} & -n_2 \frac{\rho}{\Delta_2} \\
    0 & 0 & -n_3 \frac{\rho_f}{\Delta_3} & -n_2 \frac{\rho_f}{\Delta_3} & -n_1 \frac{\rho_f}{\Delta_3} & 0 & -n_3 \frac{\rho}{\Delta_3}
  \end{pmatrix}.
\end{equation}

The matrix $\bD_v$ models the viscous dissipation, and is given by
\begin{equation} \label{eq:Dv}
  \bD_v = \begin{pmatrix}
    0 & 0 & 0 & \frac{\rho_f \eta}{\Delta_1 \kappa_1} & 0 & 0 \\
    0 & 0 & 0 & 0 & \frac{\rho_f \eta}{\Delta_2 \kappa_2} & 0 \\
    0 & 0 & 0 & 0 & 0 & \frac{\rho_f \eta}{\Delta_3 \kappa_3} \\
    0 & 0 & 0 & -\frac{\rho \eta}{\Delta_1 \kappa_1} & 0 & 0 \\
    0 & 0 & 0 & 0 & -\frac{\rho \eta}{\Delta_2 \kappa_2} & 0 \\
    0 & 0 & 0 & 0 & 0 & -\frac{\rho \eta}{\Delta_3 \kappa_3}
  \end{pmatrix}.
\end{equation}

\subsection{First-order velocity-stress system}

Combining \eqref{eq:Asv-system} and \eqref{eq:Avs-system}, and
letting the full 13-element state vector be $\bQ = \left(\begin{smallmatrix}
  \bQ_s^T & \bQ_v^T \end{smallmatrix}\right)^T$, the first-order stress-velocity
system describing three-dimensional poroelasticity is
\begin{equation} \label{eq:3d-fullsystem}
  \partial_t \bQ + \bA \bQ_{,1} + \bB \bQ_{,2} + \bC \bQ_{,3} = \bD\bQ,
\end{equation}
where
\begin{equation} \label{eq:3d-system-matrices-block}
  \begin{aligned}
  \bA &= \begin{pmatrix} 0_{7\times 7} & \bA_{sv} \\
    \bA_{vs} & 0_{6\times 6} \end{pmatrix},
  &\bB &= \begin{pmatrix} 0_{7\times 7} & \bB_{sv} \\
    \bB_{vs} & 0_{6\times 6} \end{pmatrix},\\
  \bC &= \begin{pmatrix} 0_{7\times 7} & \bC_{sv} \\
    \bC_{vs} & 0_{6\times 6} \end{pmatrix},
  &\bD &= \begin{pmatrix} 0_{7\times 7} & 0_{6\times 7} \\
    0_{7\times 6} & \bD_v \end{pmatrix}.
  \end{aligned}
\end{equation}

The reader should note that, as in~\refcartpaper{}, this system
describes low-frequency Biot poroelasticity theory --- that is, it is
valid only for angular frequencies below the critical value $\omega_c
:= \min_i \frac{\eta \phi}{\rho_f T_i \kappa_i}$.

\subsection{Energy density}

Since the constitutive relation of the poroelastic medium is linear,
the strain energy $V$
is just half the sum of the
products of the stresses with their corresponding strains,
\begin{equation} \label{eq:strain-energy-basic}
  \mathcal{V} = \frac{1}{2} \left( \sum_{I=1}^6 \tau_I e_I + p \zeta \right) =
  \frac{1}{2} \left( \btau^T \B{e} + p \zeta \right).
\end{equation}

Using equation (7.132) from~\cite{carcione:wave-book}, $\tau_I =
c_{IJ} e_J - \alpha_I p$, we can write $c_{IJ} e_J = \tau_I + \alpha_I
p$.  (Note that in the principal axes of an orthotropic material,
$\alpha_I = 0$ for $I > 3$ --- there is no equivalent shear stress in the principal axes
associated with the fluid pressure.)  Letting
$\B{S}$ be the compliance matrix of the drained skeleton --- the inverse
of the matrix formed by the drained elastic parameters
$c_{IJ}$ --- in matrix notation we have
\begin{equation} \label{eq:matrix-strain}
  \B{e} = \B{S}(\btau + p \balpha).
\end{equation}
Here $\btau$ and $\B{e}$ are arranged as $6 \times 1$ column vectors,
not as symmetric $3\times 3$ matrices.

To get the variation of fluid content $\zeta$, let us return to equation
(7.131) of~\cite{carcione:wave-book}, $p = M(\zeta - \alpha_Ie_I)$.
In matrix notation, this is $p = M(\zeta - \balpha^T\B{e}) = M(\zeta -
\balpha^T \B{S} \btau - p \balpha^T \B{S} \balpha)$.  Solving for $\zeta$ in
terms of the stress variables gives
\begin{equation} \label{eq:zeta-from-stress}
  \zeta = \left(\tfrac{1}{M} + \balpha^T \B{S} \balpha \right) p + \balpha^T \B{S} \btau.
\end{equation}

Substituting \eqref{eq:matrix-strain} and \eqref{eq:zeta-from-stress}
into \eqref{eq:strain-energy-basic}, and using the symmetry of $\B{S}$, we get
\begin{equation}
  \mathcal{V} = \tfrac{1}{2} \left( \btau^T \B{S} \btau + 2 p \balpha^T \B{S} \btau +
  p^2 \left( \tfrac{1}{M} + \balpha^T \B{S} \balpha \right) \right).
\end{equation}
In matrix form this is
\begin{equation} \label{eq:strain-energy-3d}
  \mathcal{V} = \tfrac{1}{2} \bQ_s^T \bE_s \bQ_s, \quad \bE_s := \begin{pmatrix}
    \B{S} & \B{S} \balpha \\
    \balpha^T \B{S} & \tfrac{1}{M} + \balpha^T \B{S} \balpha
  \end{pmatrix}.
\end{equation}

Meanwhile, the derivation of kinetic energy from the two-dimensional case
in~\refcartpaper{} carries over directly to three dimensions, and the
kinetic energy is
\begin{equation} \label{eq:kinetic-energy-3d}
\mathcal{T} = \tfrac{1}{2} \bQ_v^T \bE_v \bQ_v,
\end{equation}
where the matrix $\bE_v$ is
\begin{equation} \label{eq:Ev}
  \bE_v =
  \begin{pmatrix}
    \rho \B{I}_{3 \times 3} & \rho_f \B{I}_{3 \times 3} \\
    \rho_f \B{I}_{3 \times 3} & \diag(m_1, m_2, m_3)
  \end{pmatrix}.
\end{equation}

Combining \eqref{eq:strain-energy-3d} and \eqref{eq:kinetic-energy-3d}
gives the total energy per unit volume in terms of the state vector $\bQ$
as
\begin{equation} \label{eq:total-energy}
  \mathcal{E} = \mathcal{T} + \mathcal{V} = \tfrac{1}{2} \bQ^T \bE \bQ, \quad \bE := \begin{pmatrix} \bE_s & 0_{7\times 6} \\ 0_{6\times 7} & \bE_v \end{pmatrix}.
\end{equation}

We may expect $\bE$ to be positive-definite on physical grounds --- if
it were not, it would be possible to deform the medium, change its
fluid content, or set it in motion without doing work.

\subsection{Symmetrization}

In terms of the block structure of the system, $\bE \breve{\bA}$ is symmetric if
and only if $\bE_s \breve{\bA}_{sv} = (\bE_v \breve{\bA}_{vs})^T$.
After substantial algebra, it can in fact be shown that
\begin{equation}
  \bE_s \breve{\bA}_{sv} = (\bE_v \breve{\bA}_{vs})^T = - \begin{pmatrix}
    n_1 & 0 & 0 & 0 & 0 & 0 \\
    0 & n_2 & 0 & 0 & 0 & 0 \\
    0 & 0 & n_3 & 0 & 0 & 0 \\
    0 & n_3 & n_2 & 0 & 0 & 0 \\
    n_3 & 0 & n_1 & 0 & 0 & 0 \\
    n_2 & n_1 & 0 & 0 & 0 & 0 \\
    0 & 0 & 0 & -n_1 & -n_2 & -n_3
  \end{pmatrix}.
\end{equation}
Thus $\bE$ does indeed symmetrize the system.

Since $\bE$ symmetrizes the system, we immediately know that the
governing equations of three-dimensional poroelasticity are hyperbolic
by the argument of Section 2.7 of~\refcartpaper{}.  An energy norm
and energy inner product can also be defined in exactly the same
fashion as in that work.  Furthermore, we can easily see that
$\bE \bD$ is a symmetric negative-semidefinite matrix in three dimensions
as well, since
\begin{equation}
  \bE \bD =
  \begin{pmatrix}
    0_{10 \times 10} & 0_{10 \times 3} \\
    0_{3 \times 10} & -\diag \left( \eta/\kappa_1, \eta/\kappa_2, \eta/\kappa_3 \right)
  \end{pmatrix}.
\end{equation}
By the arguments of Section 2.8 of~\refcartpaper{}, almost all of the
conditions are satisfied for the energy density $\mathcal{E}$ to be a
strictly convex entropy function in the sense of Chen, Levermore, and
Liu~\cite{chen-levermore-liu:stiff-relaxation}.
There is one final condition, involving the operator $\B{\Pi}$ that
maps from the full to the reduced system; following Section 3.3
of~\refcartpaper{}, in three dimensions the matrix $\B{\Pi}$ that maps
from the full poroelastic state vector $\bQ$ to the vector of
conserved quantities of the dissipation part of the system $\B{u}$ is
\begin{equation} \label{eq:matrix-pi-3d}
  \B{\Pi} := \begin{pmatrix}
    \B{I}_{7 \times 7} & 0_{7 \times 6} \\
    0_{3 \times 7} & \B{\Pi}_v
  \end{pmatrix},
  \quad \B{\Pi}_v := \begin{pmatrix}
    1 & 0 & 0 & \rho_f/\rho & 0 & 0 \\
    0 & 1 & 0 & 0 & \rho_f/\rho & 0 \\
    0 & 0 & 1 & 0 & 0 & \rho_f/\rho
  \end{pmatrix}
\end{equation}

As in two dimensions, the fact that $\B{u}$ is conserved under the
action of the dissipation can be seen immediately because $\B{\Pi} \bD
= 0$.  The matrix $\B{G}$ that from any $\B{u}$ gives the unique
equilibrium $\bQ_{\text{eq}}$ satisfying both $\bD \bQ_{\text{eq}} =
0$ and $\B{\Pi} \bQ_{\text{eq}} = \B{u}$ is
\begin{equation}
  \bG := \begin{pmatrix}
    \B{I}_{10 \times 10} \\ 0_{3 \times 10}
  \end{pmatrix}.
\end{equation}
From these two matrices, the reduced system can be formed:
\begin{equation} \label{eq:reduced-sys-3d}
  \partial_t \B{u} + \B{\Pi} \bA \bG \B{u}_{,1} + \B{\Pi} \bB
  \bG \B{u}_{,2} + \B{\Pi} \bC \bG \B{u}_{,3} = 0.
\end{equation}

The matrix $\B{\Pi}$ allows us to see that the statements $\bD \bQ =
0$ and $\mathcal{E}'(\bQ)^T \bD \bQ = \bQ^T \bE \bD \bQ = 0$ are
equivalent to $\bQ^T \bE = \B{v}^T \B{\Pi}$ for
some $\B{v} \in \mathbb{R}^{10}$.  As in~\refcartpaper{}, if $\bQ^T
\B{E} = \B{v}^T \B{\Pi}$, we immediately have $\bQ^T \B{E} \bD \bQ =
\B{v}^T \B{\Pi} \bD \bQ = 0$ since $\B{\Pi} \bD = 0$.  Conversely, if
$\bD \bQ = 0$, the form of $\bD$ also immediately gives $q_x = q_y =
q_z = 0$.  Thus in this case the last six components of $\bQ^T \bE$
are $(\bQ^T \B{E})_8 = \rho v_x$, $(\bQ^T \B{E})_9 = \rho v_y$,
$(\bQ^T \B{E})_{10} = \rho v_z$, $(\bQ^T \B{E})_{11} = \rho_f v_x$,
$(\bQ^T \B{E})_{12} = \rho_f v_y$, and $(\bQ^T \B{E})_{13} = \rho_f
v_z$.  This implies that if $\bD \bQ = 0$, $\bQ^T \bE = \B{v}^T
\B{\Pi}$, with $\B{v}$ given by
\begin{equation}
  \B{v}^T = \begin{pmatrix}
    \bQ_s^T \bE_s & \rho v_x & \rho v_y & \rho v_z
  \end{pmatrix}.
\end{equation}
Here the stress parts of $\bQ$ and $\bE$ have been separated out for
convenience.

As in two-dimensional poroelasticity, this condition implies that
the reduced system \eqref{eq:reduced-sys-3d} is hyperbolic, and
satisfies a nonstrict subcharacteristic condition.  Equality can again
be realized in the subcharacteristic condition --- the example of this
found in Section 3.3 of~\refcartpaper{} carries over to three
dimensions --- but the argument from two dimensions that this is
harmless to the numerical solution also carries over to three dimensions.
In addition, the matrix $\B{\Pi}\breve{\bA}\bG$ again takes the form
of orthotropic elasticity, with the fluid pressure coming along as an
additional variable that does not feed back into the elastic
variables.  The flux Jacobian of the reduced system is
\begin{equation} \label{eq:reduced-sys-matrix-3d}
  \B{\Pi} \breve{\bA} \bG = \begin{pmatrix}
    0 & 0 & 0 & 0 & 0 & 0 & 0 & -n_1 c_{11}^u & -n_2 c_{12}^u & -n_3 c_{13}^u \\
    0 & 0 & 0 & 0 & 0 & 0 & 0 & -n_1 c_{12}^u & -n_2 c_{22}^u & -n_3 c_{23}^u \\
    0 & 0 & 0 & 0 & 0 & 0 & 0 & -n_1 c_{13}^u & -n_2 c_{23}^u & -n_3 c_{33}^u \\
    0 & 0 & 0 & 0 & 0 & 0 & 0 & 0 & -n_3 c_{44} & -n_2 c_{44} \\
    0 & 0 & 0 & 0 & 0 & 0 & 0 & -n_3 c_{55} & 0 & -n_1 c_{55} \\
    0 & 0 & 0 & 0 & 0 & 0 & 0 & -n_2 c_{66} & -n_1 c_{66} & 0 \\
    0 & 0 & 0 & 0 & 0 & 0 & 0 & n_1 M \alpha_1 & n_2 M \alpha_2 & n_3 M \alpha_3 \\
    -\frac{n_1}{\rho} & 0 & 0 & 0 & -\frac{n_3}{\rho} & -\frac{n_2}{\rho} & 0 & 0 & 0 & 0 \\
    0 & -\frac{n_2}{\rho} & 0 & -\frac{n_3}{\rho} & 0 & -\frac{n_1}{\rho} & 0 & 0 & 0 & 0 \\
    0 & 0 & -\frac{n_3}{\rho} & -\frac{n_2}{\rho} & -\frac{n_1}{\rho} & 0 & 0 & 0 & 0 & 0
  \end{pmatrix}.
\end{equation}

\subsection{Linear acoustics}
\label{sec:linear-acoustics}

For this work the PDEs of acoustics will be cast in the same form as
the poroelastic system \eqref{eq:3d-fullsystem}, with the same state
vector; however, in a fluid the variables $\btau$ and $\vect{v}$ will
be defined to be identically zero, as in~\refmappaper{}.  The state
variable $p$ will be used for the fluid pressure, and $\vect{q}$ for
its velocity.  The three-dimensional system's coefficient matrices
have the same block form as those for poroelasticity; the blocks are
\begin{align}
  \label{eq:A-blocks-acou-3d} \breve{\bA}_{sv} &= \begin{pmatrix}
    0_{6\times 3} & 0_{6\times 3} \\
    0_{1\times 3} & K_f \B{n}^T
  \end{pmatrix},
  &\breve{\bA}_{vs} &= \begin{pmatrix}
    0_{3\times 6} & 0_{3\times 1} \\
    0_{3\times 6} & \B{n}/\rho_f
  \end{pmatrix}.
\end{align}
The properties of a fluid are assumed isotropic, so $\breve{\bA} = n_x \bA +
n_y \bB + n_z \bC$ here is written in terms of a vector $\vect{n} =
(n_x, n_y, n_z)$ in the global problem coordinates.  (When used in
linear algebra, $\B{n}$ is taken to be a column vector.)  In a fluid, the
dissipation matrix $\bD$ is identically zero.

Just as with poroelasticity, linear acoustics also possesses an energy
density that can be expressed as a quadratic form, $\mathcal{E} =
\frac{1}{2} \bQ^T \bE \bQ$.  The energy divides neatly into kinetic
and potential, giving the same block structure for $\bE$ as before,
with the blocks for acoustics equal to
\begin{equation} \label{eq:E-blocks-acou-3d}
  \bE_s = \begin{pmatrix}
    0_{6 \times 6} & 0_{6 \times 1} \\
    0_{1 \times 6} & 1/K_f
  \end{pmatrix}, \quad
  \bE_v = \begin{pmatrix}
    0_{3 \times 3} & 0_{3 \times 3} \\
    0_{3 \times 3} & \rho_f \B{I}_{3 \times 3}
  \end{pmatrix}
\end{equation}
This matrix $\bE$ is only positive-semidefinite, not positive-definite
as for poroelasticity.  Similarly to the two-dimensional case
discussed in~\refmappaper{}, however, the null space of $\bE$ consists
only of the variables that are defined to be identically zero in the
fluid.  This means $\bE$ is essentially positive-definite, and can
still be used to define an energy inner product and norm for
acoustics.  Just as for poroelasticity, $\bE$ symmetrizes the
first-order hyperbolic system for acoustics.

\subsection{Interface conditions}
\label{sec:interface-cond}

The same interface conditions are used here in three dimensions as
were used in two dimensions in~\refmappaper{} --- and in fact the same vector formulas may be
used, since the ones written in Section 2.4 of~\refmappaper{} are
equally valid in two or three dimensions.

To reprise, open pores,
closed pores, or imperfect hydraulic contact between two poroelastic
media may be expressed as the set of conditions
\begin{equation} \label{eq:poro-poro-physical}
  \begin{aligned}
    \btau_l \cdot \vect{n} &= \btau_r \cdot \vect{n}\\
    \vect{v}_l &= \vect{v}_r\\
    \vect{q}_l \cdot \vect{n} &= \vect{q}_r \cdot \vect{n}\\
    \eta_d (p_l - p_r) &= Z_f (1 - \eta_d) \widehat{\vect{q} \cdot \vect{n}},
  \end{aligned}
\end{equation}
where the subscripts $l$ and $r$ denote the arbitrarily-chosen left
and right sides of the interface, the vector $\vect{n}$ is the unit
interface normal pointing from the left medium to the right one, $Z_f$
is the acoustic impedance of the fluid in the left medium, and
$\eta_d$ is the {\em interface discharge efficiency}, a nondimensional
measure of the resistance of the interface to fluid flow.  In this
formulation, $\eta_d = 1$ corresponds to fluid flowing across the
interface unhindered, while $\eta_d = 0$ corresponds to a completely
impermeable interface, and values of $\eta_d$ between 0 and 1 indicate
an interface that allows fluid to pass, but only if driven by a
pressure difference.
The quantity $\widehat{\vect{q} \cdot \vect{n}}$ is equal to both
$\vect{q}_l \cdot \vect{n}$ and $\vect{q}_r \cdot \vect{n}$.
Similar interface conditions between a poroelastic medium and a fluid
may be expressed as
\begin{equation} \label{eq:poro-fluid-physical}
  \begin{aligned}
    \vect{q}_f \cdot \vect{n} &= (\vect{v}_p + \vect{q}_p)
    \cdot \vect{n}\\
    -p_f \vect{n} &= \btau_p \cdot \vect{n}\\
    \eta_d (p_p - p_f) &=  Z_f (1 - \eta_d) \vect{q}_p \cdot \vect{n},
  \end{aligned}
\end{equation}
where $Z_f$ is the impedance of the fluid medium, the subscripts $p$
and $f$ denote the poroelastic and fluid media, and the
unit vector $\vect{n}$ points from the poroelastic medium into the fluid.

\section{Finite volume methods for mapped grids in three dimensions}
\label{sec:mapped-fvm-3d}

Despite the increase in the number of spatial dimensions and the
expansion of the state vector from 8 to 13 elements, most of the
details of the numerical method for three-dimensional poroelasticity
and \porofluid{} systems on mapped grids are closely analogous to the
two-dimensional methods of~\refcartpaper{} and~\refmappaper{}.
It is primarily the changes between the
two-dimensional and three dimensional method that are described here.
Section \ref{sec:mapped-grids-3d} describes how mapped grids are used
in three dimensions, while in Section \ref{sec:riemann-solve-3d} the
Riemann solution process is discussed, including the interface
condition matrices used to solve the Riemann problem between a fluid
and a poroelastic medium in three dimensions.

There are also some algorithmic details that differ between two
dimensions and three.  Specifically, the combination of mapped grids
and anisotropic materials requires a change to the limiter algorithm
in order to cope with potential difficulties in consistently tracking
the polarization of the
shear waves; this change is discussed and evaluated in two dimensions
in Section \ref{sec:elimiter}.  Also, since the transverse Riemann
solutions for three-dimensional high-resolution finite volume methods
are both computationally expensive and time-consuming to program,
dimensional splitting becomes very attractive.  While dimensional
splitting is by no means a new algorithmic development, Section
\ref{sec:dimsplit} gives a simple overview.  Finally, Section
\ref{sec:software-3d} gives a brief summary of the software frameworks
in which these algorithms are implemented.

\subsection{Mapped grids in three dimensions}
\label{sec:mapped-grids-3d}

As in~\refmappaper{}, the three-dimensional numerical solution
procedure uses logically rectangular mapped grids.  In three
dimensions, though, mapped grid quantities such as interface normals
and cell face areas are more difficult to define.

Each cell in the mapped grid is defined by the physical coordinates of
its vertices, which are computed by applying the mapping function to
the vertices of the cell in the computational domain.  Physical
coordinates will be denoted by $x$, $y$, and $z$, or the position
vector $\vect{r}$, while computational
coordinates are $\xi_1$, $\xi_2$, and $\xi_3$.  From its
vertices, cells are defined by a trilinear mapping.  Let
$\eta_1$, $\eta_2$, and $\eta_3$ be cell-local computational
coordinates, defined from the global computational coordinates by
\begin{equation}
  \eta_i = \frac{\xi_i - \xi_{i0}}{\Delta \xi_i},
\end{equation}
where $\xi_{i0}$ is the lowest extent of the cell in global
computational coordinate $i$ and $\Delta \xi_i$ is the grid spacing in
computational coordinate $i$.  In local computational coordinates, the
cell is thus the unit cube $[0,1]^3$.  From these coordinates, the
cell is parameterized in physical coordinates by
\begin{equation} \label{eq:3d-parametric-cell}
  \vect{r}(\eta_1, \eta_2, \eta_3) = \sum_{i=1}^2 \sum_{j=1}^2
  \sum_{k=1}^2 \vect{r}_{ijk} N_i(\eta_1) N_j(\eta_2) N_k(\eta_3),
\end{equation}
where the $\vect{r}_{ijk}$ are the vertices, with each subscript
denoting the position in the corresponding computational direction (so
for example $\vect{r}_{112}$ is the $-\xi_1$, $-\xi_2$, $+\xi_3$
vertex), and the $N$ functions are
$N_1(\eta) = 1-\eta$, $N_2(\eta) = \eta$.

Defining the normal vector to a cell face stretched between four
essentially arbitrary vertices is not trivial, because the vertices
may not be coplanar.  A sensible requirement, though, seems to be
for the normal $\vect{n}$ and area $A$ of a face to
satisfy
\begin{equation} \label{eq:3d-normal-condition}
  \vect{n} A = \int_{\text{face}} \vect{n}_{\text{local}}\,dA,
\end{equation}
where $\vect{n}_{\text{local}}$ is the local unit normal at each point
on the face.  If the normals point outward, this implies
\begin{equation}
  \sum_{i\, \in\, \text{faces}} \vect{n}_i A_i = \int_{\partial
    (\text{cell})} \vect{n}_{\text{local}}\,dA = 0,
\end{equation}
a fundamental property of a closed surface.  In particular, for a
conservative finite volume method this implies that the net flux of a
constant vector through the cell is zero, which it should be
since the divergence of a constant vector is zero.

The simplest way to satisfy \eqref{eq:3d-normal-condition} is to
calculate the integral on the right, then let $A$ be the magnitude of
the resulting vector and $\vect{n}$ be the unit vector in that
direction.  For the parameterization \eqref{eq:3d-parametric-cell},
this integral reduces to the cross product of the vectors connecting
the midpoints of the face edges, so, for example, on the $+\xi_1$ face
we get
\begin{equation}
  \vect{n} A = \left( \tfrac{1}{2} \left( \vect{r}_{221} +
  \vect{r}_{222} \right) - \tfrac{1}{2} \left( \vect{r}_{211} +
  \vect{r}_{212} \right) \right) \times
  \left( \tfrac{1}{2} \left( \vect{r}_{212} + \vect{r}_{222} \right) -
  \tfrac{1}{2} \left( \vect{r}_{211} + \vect{r}_{221} \right) \right).
\end{equation}

The cell volume, needed for calculation of the capacity $\kappa$, is
more difficult to compute.  While an analytical expression can be
developed for the integral of the Jacobian of
\eqref{eq:3d-parametric-cell}, this expression is quite cumbersome,
and it is easier to evaluate the integral by quadrature.  The Jacobian
is quadratic in each local coordinate $\eta_i$, so a tensor product of
two-point Gauss-Legendre quadrature in each direction
evaluates it to machine precision.  The cell centroid location,
which is not a fundamental part of the finite volume method but is
useful for evaluating spatially-varying initial and boundary
conditions, can be calculated the same way --- since the position
vector \eqref{eq:3d-parametric-cell} is first-order in each local
coordinate, its product with the Jacobian is at most third-order, so
two-point Gauss-Legendre quadrature still evaluates it
exactly.

\subsection{Riemann problems on three-dimensional mapped grids}
\label{sec:riemann-solve-3d}

The solution process for Riemann problems on three-dimensional mapped
grids is very similar to the process on two-dimensional mapped grids
detailed in~\refmappaper{}.  This section will
primarily focus on the changes necessary in passing to three
dimensions.

\subsubsection{Eigenvalues and eigenvectors}
\label{sec:poro-eigen-3d}

As in two dimensions, the eigenvectors for three-dimensional acoustics
are quite simple; for the matrix $\breve{\bA} = n_x \bA + n_y \bB +
n_z \bC$ the vectors for left- and right-going waves may easily be verified as
\begin{equation}
  \begin{aligned}
    \B{r}_{\text{acoustic, left}} &= \begin{pmatrix}
      0 & 0 & 0 & 0 & 0 & 0 & -Z_f & 0 & 0 & 0 & n_x & n_y & n_z
    \end{pmatrix}^T\\
    \B{r}_{\text{acoustic, right}} &= \begin{pmatrix}
      0 & 0 & 0 & 0 & 0 & 0 & Z_f & 0 & 0 & 0 & n_x & n_y & n_z
    \end{pmatrix}^T.
  \end{aligned}
\end{equation}

For three-dimensional poroelasticity, the eigenvectors may be found by
a procedure very similar to that used in two dimensions
in~\refmappaper{}, converting the problem to a symmetric generalized
eigenproblem and exploiting the block structure of $\breve{\bA}$.  A
detailed account of this process is given in Section 8.3.1
of~\cite{lemoine:thesis}.

\subsubsection{Solution of the Riemann problem}
\label{sec:interface-3d}

For the case of identical materials on either side of an interface,
$\bE$-orthogonality of the eigenvectors allows easy extraction of the
wave strengths from the difference in states, just as
in~\refmappaper{}.  For the case of different materials, with an
interface condition between them, the overall solution procedure is
the same as in that paper, but some of the specifics differ because
of the different state vector.

Because there is nothing specific to a having particular number of
spatial dimensions in the solution procedure for Riemann problems with
interface conditions in~\refmappaper{}, the same overall procedure can
be used here, the only differences being the size of the state vector
$\bQ$ and the number of waves.  The only task remaining is to write
the matrices $\bC_l$ and $\bC_r$ corresponding to the interface
conditions \eqref{eq:poro-poro-physical} and
\eqref{eq:poro-fluid-physical}.

The fluid-poroelastic interface condition will be treated first.
Taking the left medium to be poroelastic, and introducing the
parameter $Z' := Z_f (1-\eta_d)$ for brevity, a component-by-component
accounting of the correspondence between physical variables and the
entries of $\bQ$ gives
\begin{multline}
    \bC_{l,\text{poro-fluid}} \\= \begin{pmatrix}
      0 & 0 & 0 & 0 & 0 & 0 & 0 &    
      n_x & n_y & n_z & n_x & n_y & n_z \\    
      n_x & 0 & 0 & 0 & n_z & n_y & 0 &    
      0 & 0 & 0 & 0 & 0 & 0 \\    
      0 & n_y & 0 & n_z & 0 & n_x & 0 &    
      0 & 0 & 0 & 0 & 0 & 0 \\    
      0 & 0 & n_z & n_y & n_x & 0 & 0 &    
      0 & 0 & 0 & 0 & 0 & 0 \\    
      0 & 0 & 0 & 0 & 0 & 0 & \eta_d &    
      0 & 0 & 0 & -Z' n_x & -Z' n_y & -Z' n_z
    \end{pmatrix}
\end{multline}
\begin{equation}
    \bC_{r,\text{poro-fluid}} = \begin{pmatrix}
      0 & 0 & 0 & 0 & 0 & 0 & 0 &    
      0 & 0 & 0 & n_x & n_y & n_z \\    
      0 & 0 & 0 & 0 & 0 & 0 & -n_x &    
      0 & 0 & 0 & 0 & 0 & 0 \\    
      0 & 0 & 0 & 0 & 0 & 0 & -n_y &    
      0 & 0 & 0 & 0 & 0 & 0 \\    
      0 & 0 & 0 & 0 & 0 & 0 & -n_z &    
      0 & 0 & 0 & 0 & 0 & 0 \\    
      0 & 0 & 0 & 0 & 0 & 0 & \eta_d &    
      0 & 0 & 0 & 0 & 0 & 0
    \end{pmatrix}.
\end{equation}
The vector $\B{n} = (n_x, n_y, n_z)$ is the unit
interface normal pointing from the poroelastic medium into the fluid;
$Z_f$ is the fluid acoustic impedance, and $\eta_d$ is the interface
discharge efficiency, defined in Section \ref{sec:interface-cond}.  If
the poroelastic medium is on the right, the
subscripts $l$ and $r$ may simply be exchanged and the normal $\B{n}$
negated.

For the poroelastic-to-poroelastic interface condition, as
in~\refmappaper{}, write the quantity $\widehat{\vect{q} \cdot
  \vect{n}}$ from \eqref{eq:poro-poro-physical} as a
weighted average of the normal flow rates on both sides of the
interface, $\widehat{\vect{q} \cdot \vect{n}} = (1-\zeta) \vect{q}_l
\cdot \vect{n} + \zeta \vect{q}_r \cdot \vect{n}$.  Then the interface
condition matrices $\bC_l$ and $\bC_r$ become
\begin{multline}
    \bC_{l,\text{poro-poro}} \\= \begin{pmatrix}
      n_x & 0 & 0 & 0 & n_z & n_y & 0 &    
      0 & 0 & 0 & 0 & 0 & 0 \\    
      0 & n_y & 0 & n_z & 0 & n_x & 0 &    
      0 & 0 & 0 & 0 & 0 & 0 \\    
      0 & 0 & n_z & n_y & n_x & 0 & 0 &    
      0 & 0 & 0 & 0 & 0 & 0 \\    
      0 & 0 & 0 & 0 & 0 & 0 & 0 &    
      1 & 0 & 0 & 0 & 0 & 0 \\    
      0 & 0 & 0 & 0 & 0 & 0 & 0 &    
      0 & 1 & 0 & 0 & 0 & 0 \\    
      0 & 0 & 0 & 0 & 0 & 0 & 0 &    
      0 & 0 & 1 & 0 & 0 & 0 \\    
      0 & 0 & 0 & 0 & 0 & 0 & 0 &    
      0 & 0 & 0 & n_x & n_y & n_z \\    
      0 & 0 & 0 & 0 & 0 & 0 & \eta_d &    
      0 & 0 & 0 & -Z_l' n_x & -Z_l' n_y & -Z_l' n_z    
    \end{pmatrix}
\end{multline}
\begin{equation}
    \bC_{r,\text{poro-poro}} = \begin{pmatrix}
      n_x & 0 & 0 & 0 & n_z & n_y & 0 &    
      0 & 0 & 0 & 0 & 0 & 0 \\    
      0 & n_y & 0 & n_z & 0 & n_x & 0 &    
      0 & 0 & 0 & 0 & 0 & 0 \\    
      0 & 0 & n_z & n_y & n_x & 0 & 0 &    
      0 & 0 & 0 & 0 & 0 & 0 \\    
      0 & 0 & 0 & 0 & 0 & 0 & 0 &    
      1 & 0 & 0 & 0 & 0 & 0 \\    
      0 & 0 & 0 & 0 & 0 & 0 & 0 &    
      0 & 1 & 0 & 0 & 0 & 0 \\    
      0 & 0 & 0 & 0 & 0 & 0 & 0 &    
      0 & 0 & 1 & 0 & 0 & 0 \\    
      0 & 0 & 0 & 0 & 0 & 0 & 0 &    
      0 & 0 & 0 & n_x & n_y & n_z \\    
      0 & 0 & 0 & 0 & 0 & 0 & \eta_d &    
      0 & 0 & 0 & Z_r' n_x & Z_r' n_y & Z_r' n_z    
    \end{pmatrix}.
\end{equation}
As with the fluid-poroelastic matrices, the parameters $Z'_l :=
(1-\zeta) Z_f (1-\eta_d)$ and $Z_r' := \zeta Z_f (1-\eta_d)$ have been
introduced for convenience.  Based on the results of Section
3.2 of~\refmappaper{}, $\zeta = \frac{1}{2}$ is used in all
cases.

Aside from these new interface condition matrices, the Riemann solution
procedure in three dimensions is identical to that in two dimensions
in~\refmappaper{}.

\subsection{Shear waves and revised limiter algorithm}
\label{sec:elimiter}

One subtle difficulty in using high-resolution finite volume methods
for three dimensional elasticity and poroelasticity comes with the
shear waves.  In three dimensions, for any given propagation direction
there are two possible polarizations of shear wave, which in an orthotropic material may or may
not have the same speed.  The possibility
of different speeds obliges a general-purpose
code to treat them as two distinct waves, but if the
speeds are the same, the matrix $\breve{\bA}$ has a two-dimensional
eigenspace, without any intrinsic way to assign waves to one
family or the other.  This still does not present difficulties in the
core Lax-Wendroff method, but applying a limiter to a wave requires
comparing it to the upwind wave in the same family.  On a rectilinear
grid it would be possible to arbitrarily define (say)
vertically-polarized shear waves to be one family and
horizontally-polarized ones to be the other, but on a general mapped
grid it is impossible to choose polarization directions in a way that
smoothly varies over all possible cell interface directions --- the
popular result that ``you can't comb the hair on a sphere'' --- so
there could be discontinuities in the chosen polarization directions
from one Riemann problem to the next.  The limiter would see these as
solution discontinuities, and would act to suppress the higher-order
terms of the method around them, even if the solution were in fact
smooth.  This combination of possibly anisotropic materials and mapped
grids makes it challenging to formulate a good wave limiting
algorithm.

The solution used here is to exploit the $\bE$-orthogonality of the
eigenvectors of $\breve{\bA}$ to find the component of the upstream
waves in the direction of the wave to be limited, rather than using
the wave family number.  In the classical approach to wave limiting,
the wave strength ratio $\theta$ for wave $m$ at cell interface
$(i-1/2, j, k)$ is computed as
\begin{equation} \label{eq:3dlim-theta-classical}
  \theta_{\text{classical}} :=
  \frac{(\mathcal{W}^m_{i-1/2,j,k})^T \mathcal{W}^m_{I-1/2,j,k}}
  {(\mathcal{W}^m_{i-1/2,j,k})^T \mathcal{W}^m_{i-1/2,j,k}},
\end{equation}
where interface $(I-1/2,j,k)$ is the upwind interface.  In its most
minimal form, applied only to the shear waves, the new {\em energy
  inner product limiter} ($\bE$-limiter for short) replaces the
unweighted inner products with energy inner products, and takes the
inner product of the shear wave to be limited with the sum of the
upwind shear waves; if the shear waves are identified as S1 and
S2, the wave strength ratio is
\begin{equation} \label{eq:3dlim-theta-shear}
  \theta_{\bE,\text{shear}} := \frac{(\mathcal{W}^m_{i-1/2,j,k})^T
    \bE (\mathcal{W}^{\text{S1}}_{I-1/2,j,k}
    + \mathcal{W}^{\text{S2}}_{I-1/2,j,k})}
  {(\mathcal{W}^m_{i-1/2,j,k})^T \bE \mathcal{W}^m_{i-1/2,j,k}}.
\end{equation}
For an inhomogeneous domain, $\bE$ is the energy density matrix for
the medium into which the wave to be limited is propagating.  Because
the wave eigenvectors are $\bE$-orthogonal, the numerator of
\eqref{eq:3dlim-theta-shear} gives exactly the
component of the upwind shear waves in the direction of the wave being
limited, regardless of the choice of eigenvectors and assignment of
wave family numbers at the upwind interface.  Once the wave strength
ratio $\theta$ has been calculated, the limiter function
$\phi(\theta)$ is applied and the wave is scaled accordingly, exactly
as in~\cite{rjl:fvm-book}.

There is a potential implementation difficulty with the $\bE$-limiter
wave strength as computed in \eqref{eq:3dlim-theta-shear}.  As
written, the formula for computing $\theta$ requires knowing which
upstream waves are the shear waves.  Normally the shear waves are
intermediate in speed between the fast and slow P waves, but that does
not necessarily have to be the case.  If the shear modulus is
extremely low, for example, both shear waves might be slower than
the slow P wave.  While this difficulty could be avoided by
formulating a procedure that explicitly computes the shear and
longitudinal wave eigenvectors separately, such a procedure would be
complex to implement in the context of waves propagating in an
arbitrary direction through an orthotropic medium, since the shear and
extensional deformation are coupled in any set of axes other than the
material principal axes.  Rather than attempt to determine which wave is
which, here the $\bE$-orthogonality is again exploited by adding
together all the upstream waves.  The full
form of the $\bE$-limiter wave strength ratio as implemented in the
three-dimensional simulation code is then
\begin{equation} \label{eq:3dlim-theta}
  \theta_{\bE} := \frac{(\mathcal{W}^m_{i-1/2,j,k})^T \bE
    \sum_{n, \text{same direction}} \mathcal{W}^n_{I-1/2,j,k}}
  {(\mathcal{W}^m_{i-1/2,j,k})^T \bE
    \sum_{n, \text{same direction}} \mathcal{W}^n_{i-1/2,j,k}},
\end{equation}
where the sums are over waves going in the same direction (left or
right) as the wave to be limited.  The expression in the denominator
is a more efficient way of calculating $(\mathcal{W}^m_{i-1/2,j,k})^T
\bE \mathcal{W}^m_{i-1/2,j,k}$ for multiple waves --- the two
expressions are equal by $\bE$-orthogonality of the waves, but
computing $\bE \sum \mathcal{W}^n_{i-1/2,j,k}$, storing it, and successively
computing its inner product with each $\mathcal{W}^m_{i-1/2,j,k}$
requires fewer floating-point operations than computing
$(\mathcal{W}^m_{i-1/2,j,k})^T \bE \mathcal{W}^m_{i-1/2,j,k}$ for each
wave separately.

The $\bE$-orthogonality of the waves means that the new formula
\eqref{eq:3dlim-theta} gives the
same results as the original formula \eqref{eq:3dlim-theta-classical}
for fast and slow P waves for homogeneous media if
successive cell interfaces are parallel.  However, because the
$\breve{\bA}$ matrix depends on the normal direction of the grid
interface, on more general mapped grids a wave in one family may not be
$\bE$-orthogonal to waves in different families at the upwind
interface.  It may in fact be appropriate to include contributions
from other wave families in such cases --- for instance, a plane wave
in one direction in a single wave family may have components in multiple families
if expressed in terms of eigenvectors of an $\breve{\bA}$ matrix
computed for a different direction.  In
order to assess the effect of this, the
cylindrical scatterer test cases of~\refmappaper{} are re-run here
with the original limiter formulation and the $\bE$-limiter, both
using the Monotonized Centered (MC) limiter function for $\phi(\theta)$.  All 18 scatterer
cases are examined.  Because the $\bE$-limiter adds together waves
from different families, it seems appropriate to also test it in
combination with the $f$-wave approach of Bale, LeVeque, Mitran, and
Rossmanith~\cite{bale-rjl-mitran-rossmanith:f-waves}, to see whether
weighting the different waves by their speeds in the sum would give
noticeably different results.  The original limiter formulation is not
run here with $f$-waves because the poroelastic material in the scatterer model
is isotropic --- the $f$-wave formulation weights waves by their
speeds, and because the wave speeds are the same in all directions,
this weighting would have no effect on the wave strength ratios when
comparing within the same wave family.

\begin{table}
  \caption[Effect of using the energy inner product
    limiter]{Percentage change in error caused by using the
    $\bE$-limiter, with ordinary waves or $f$-waves, on the
    cylindrical scatterer test cases of~\refmappaper{}.  Percentage
    change is relative to error in solution obtained using the
    conventional wave strength ratio \eqref{eq:3dlim-theta-classical}.
    Statistics are computed across all 18 scatterer test cases.}
  \label{tab:elimiter-sc-results}
  \begin{center}
    \begin{tabular}{lrrrr}
      \toprule
      & \multicolumn{2}{c}{Ordinary waves} & \multicolumn{2}{c}{$f$-waves} \\
      \cmidrule(r){2-3} \cmidrule(r){4-5}
      & 1-norm & Max-norm & 1-norm & Max-norm \\
      \midrule
      \multicolumn{5}{l}{All grids}\\
      \midrule
      Maximum & $+3.42$\% & $+2.86$\% & $+61.95$\% & $+10.24$\%\\
      Minimum & $-8.21$\% & $-21.43$\% & $-2.48$\% & $-22.48$\%\\
      Mean & $-1.97$\% & $-1.54$\% & $+5.77$\% & $-0.28$\%\\
      Median & $-2.05$\% & $-0.62$\% & $+2.40$\% & $-0.22$\%\\
      \midrule
      \multicolumn{5}{l}{Finest grid}\\
      \midrule
      Maximum & $+3.21$\% & $+1.60$\% & $+61.95$\% & $+5.41$\%\\
      Minimum & $-7.75$\% & $-8.55$\% & $-2.48$\% & $-9.49$\%\\
      Mean & $-3.05$\% & $-0.63$\% & $+9.55$\% & $-0.35$\%\\
      Median & $-3.65$\% & $+0.23$\% & $+5.63$\% & $+1.23$\%\\
      \bottomrule
    \end{tabular}
  \end{center}
\end{table}

Table \ref{tab:elimiter-sc-results} lists the percent changes in error
due to incorporating the $\bE$-limiter with or without $f$-waves.
Adding the $\bE$-limiter by itself typically gives a modest reduction
in error, with some cases seeing substantial reduction and others
seeing slight increases.  Using the $\bE$-limiter in combination with
$f$-waves, on the other hand, typically increases the error, sometimes
dramatically.  Because of this, the 3D simulation code uses the
$\bE$-limiter wave strength ratio \eqref{eq:3dlim-theta} with ordinary
waves, not $f$-waves.

As a final comment, note that while the new wave strength ratio
\eqref{eq:3dlim-theta} is formulated using additional knowledge of the
structure of the poroelasticity system compared to the classical
strength ratio \eqref{eq:3dlim-theta-classical} --- namely knowledge
of the $\bE$-orthogonality of the waves --- and thus is not as easy to
generalize, similar formulas could be constructed for other systems if
similar orthogonality relations can be found.  This would be
advantageous for applying high-resolution finite volume methods to
hyperbolic systems where it is not always clear which wave should be
compared to which when applying the limiter.  The reader should also
note that, because this is a new way to calculate the wave strength
ratio $\theta$, it is compatible with any limiter function
$\phi(\theta)$, and can be used in conjunction with recent advances in
limiter functions such as those reviewed by
Kemm~\cite{kemm:limiters}.

\subsection{Dimensional splitting}
\label{sec:dimsplit}

Transverse wave propagation in three dimensions is both complex and
computationally expensive.  Beyond the normal Riemann solve, which is
always necessary, the classical 3D transverse propagation approach
worked out by Langseth and LeVeque~\cite{langseth-rjl:3d-transverse}
and implemented in \clawpack{} requires eight transverse Riemann
solves per cell interface, and in addition eight double-transverse
Riemann solves, which provide third-order terms that are necessary for
stability.  Extending the new two-dimensional transverse solve scheme
of~\refmappaper{} to three dimensions in an
analogous fashion would require 16 transverse Riemann solves and
perhaps as many as 32 double-transverse solves; including the normal
solve, this could be as many as 49 Riemann solves per interface, a
prohibitively high computational cost.  Most of the computational
effort in the poroelastic Riemann solver is in a lengthy setup phase
where eigenvectors and coefficient matrices are computed, and since
this phase does not depend on the cell states or fluctuations, it
would be possible to bundle together all of the transverse or
double-transverse solves for a particular interface, and solve them
all together for a cost only marginally higher than a single solve.
This would reduce the number of times the setup phase is run to three
times per cell interface, but it would require a substantial rewrite
of the \clawpack{} internals, which would be prohibitively
time-consuming and error-prone.

Because of the computational expense of the transverse solves, all
three\hyp dimensional results in this work are run using dimensional
splitting.  For the dimensionally-split approach the normal Riemann
problems are solved in only one grid direction at a time.  Their solutions are
used to update the cells to an intermediate state, and this
intermediate state is used to solve the normal Riemann problems in the
next grid direction; the results are used to update the cells to a new
intermediate state, which is used to solve the normal Riemann problems
in the final grid direction and update the cells to the next time
step.  In combination with Strang splitting for the source term, then,
the procedure to advance the solution by $\Delta t$ from time step $n$
to $n+1$ with dimensional splitting runs as follows:
\begin{enumerate}
\item Advance $\bQ^n$ by $\Delta t/2$ using the source
  term, giving $\bQ^{(0)}$.
\item Advance $\bQ^{(0)}$ by $\Delta t$ using Riemann solves in the
  $i$ direction, giving $\bQ^{(1)}$.
\item Advance $\bQ^{(1)}$ by $\Delta t$ using Riemann solves in the
  $j$ direction, giving $\bQ^{(2)}$.
\item Advance $\bQ^{(2)}$ by $\Delta t$ using Riemann solves in the
  $k$ direction, giving $\bQ^{(3)}$.
\item Advance $\bQ^{(3)}$ by $\Delta t/2$ using the source term
  again.  The result is $\bQ^{n+1}$.
\end{enumerate}
While it is only first-order accurate, this work nonetheless
uses dimensional splitting exclusively, because it appears to be the
only timely way to obtain numerical solutions for these cases, both in
terms of software development time and program execution time.

\subsection{Numerical software}
\label{sec:software-3d}

The numerical solution techniques described here were implemented
using a hybrid of several different versions of the \clawpack{} finite
volume software.  A pure-Fortran implementation was written to
interface with \clawpack{} 4.3, which was the last version before
\clawpack{} 5.0 (which was not yet available as of this writing) to
support three-dimensional problems.  A hybrid Python-Fortran
implementation was also written for \textsc{PyClaw}~\cite{pyclaw} in
order to be able to use the \textsc{PetClaw}~\cite{petclaw} variant of
\textsc{PyClaw} to run in parallel on large workstation-class
computers or clusters.

\section{Results}
\label{sec:results}

With the numerical methods formulated for three-dimensional
poroelasticity and \porofluid{} systems, it is now time to apply these
methods to some test problems.  Section \ref{sec:plane-soln-3d}
details the construction of plane wave solutions analogous to those of
Section 4.1 of~\refcartpaper{}.  Section \ref{sec:plane-conv-3d} then
uses these solutions to examine the convergence behavior of the
numerical methods of Section \ref{sec:mapped-fvm-3d}, and section
\ref{sec:elimiter-test} examines the performance the new
$\bE$-limiter.  Following this, the results of a demonstration problem
that exercises almost all of the functionality of the
three-dimensional code are presented in section \ref{sec:3d-demo}.

\subsection{Analytic plane wave solution}
\label{sec:plane-soln-3d}

The process of generating a plane wave solution begins by prescribing
a unit vector $\Bell$ in the direction of the desired wavevector, an
angular frequency $\omega$, the orientation of the principal axes of
the material, a desired wave family, and for shear waves a desired
polarization direction $\vect{s}$.  Given these inputs, first the
vectors $\Bell$ and $\vect{s}$ are transformed into the material
principal coordinates to simplify subsequent calculations.  Following
this, the complex wavenumber $k$ and wave eigenvector
$\B{v}$ are obtained by taking an ansatz for the solution $\bQ$ of the
form
\begin{equation} \label{eq:3d-plane-ansatz}
  \bQ = \B{v} \exp \left( i ( k (\ell_1 x_1 + \ell_2 x_2 + \ell_3 x_3) - \omega
  t) \right).
\end{equation}
Here $\ell_1$, $\ell_2$, and $\ell_3$ are the components of $\Bell$ in the
material principal coordinates, and $x_1$, $x_2$, and $x_3$ are
distances along the principal material axes.

Substituting this ansatz into
the first-order system for three-dimensional poroelasticity
\eqref{eq:3d-fullsystem} results in the eigenproblem
\begin{equation} \label{eq:3d-plane-eigen-base}
  -i \omega \B{v} + ik \breve{\bA} \B{v} = \bD \B{v},
\end{equation}
where as usual $\breve{\bA} = \ell_1 \bA + \ell_2 \bB + \ell_3 \bC$.
Rearranging and rescaling by making the substitution $\B{v} =
\bE^{-1/2} \B{w}$, then multiplying from the left by $\bE^{-1/2}$ results
in the complex symmetric generalized eigenproblem
\begin{equation} \label{eq:3d-plane-eigen}
  \bE^{1/2} \breve{\bA} \bE^{-1/2} \B{w} = k^{-1} \left( \omega \B{I} - i
  \bE^{1/2} \bD \bE^{-1/2} \right) \B{w}.
\end{equation}
This second form of the eigenproblem is easier to work with
numerically --- multiplying by the square root and inverse square root
of the energy density matrix $\bE$ improves the relative scaling of
the components of $\breve{\bA}$, and letting $k^{-1}$ be the
eigenvalue allows the null vectors of $\breve{\bA}$ to correspond to
zero eigenvalues rather than the infinite ones that would result if
$k$ were placed in the role of eigenvalue.

From the solutions of the eigenproblem \eqref{eq:3d-plane-eigen}, the
eigenvector $\B{w}$ and complex wavenumber $k$ are extracted that
correspond to the desired wave family.  The wave eigenvector $\B{v}$
is then computed using $\B{v} = \bE^{-1/2} \B{w}$.  In the case where
a shear wave is requested and the two shear wave speeds are equal, the
vector from the two-dimensional shear eigenspace is chosen that has
solid velocity as close to parallel to the prescribed 
direction $\vect{s}$ as possible.  The vector $\B{v}$ is then
normalized to unit $\bE$-norm, and its complex phase is adjusted so
that the dot product of its solid velocity component with a reference
direction --- $\Bell$ for fast and slow P waves, $\vect{s}$ for shear
waves --- is pure real and positive.  Finally, since the eigenproblem
was solved using the system matrices for the principal material axes,
$\vect{v}$ is transformed back into the global computational axes.

\subsection{Plane wave convergence studies}
\label{sec:plane-conv-3d}

As in~\refcartpaper{} in two dimensions, the three-dimensional code
is first tested using analytical plane wave solutions.  Based
on the results of~\refcartpaper{}, and because of the high
computational cost of three-dimensional simulation, only viscous
high-frequency test cases are run here.  Good convergence behavior for
these cases implies that the underlying wave propagation algorithm
would also perform well for inviscid cases, and from~\refcartpaper{}
we already know to expect first-order convergence for low-frequency
viscous cases, regardless of how well the code would perform
otherwise.  Even with the restriction to viscous high-frequency cases,
only a relatively small number of cases are examined in order to keep
the computational cost of these convergence studies reasonable.

\begin{table}
\caption{Properties of the orthotropic sandstone used in test cases,
  take from de la Puente et
  al.~\cite{delapuente-dumbser-kaser-igel:poro-dg}.  Wave speeds are
  correct in the high-frequency limit; $c_{pf}$ is the fast P wave
  speed, $c_s$ is the S wave speed, $c_{ps}$ is the slow P wave speed,
  and $\tau_d$ is the time constant for dissipation.  Subscript
  numbers indicate principal directions.  The material
  is isotropic in the 1-2 plane, so the 2-direction properties are
  related to the tabulated values by $c_{22} = c_{11}$, $c_{23} =
  c_{13}$, $c_{44} = c_{55}$, $c_{66} = (c_{11} - c_{12})/2$,
  $\kappa_2 = \kappa_1$, and $T_2 = T_1$.}
\label{tab:matprops}
\begin{center}
  \begin{tabular}{rrrrrr}
    \toprule
    \multicolumn{4}{c}{Base properties} & \multicolumn{2}{c}{Derived properties} \\
    \cmidrule(r){1-4} \cmidrule(r){5-6}
    $K_s$ & 80 GPa & $\kappa_1$ & $600 \times 10^{-15}$ m$^2$ & $c_{pf1}$ & 6000 m/s \\
    $\rho_s$ & 2500 kg/m$^3$ & $\kappa_3$ & $100 \times 10^{-15}$ m$^2$ & $c_{pf3}$ & 5260 m/s \\
    $c_{11}$ & 71.8 GPa & $T_1$ & 2 & $c_{s1}$ & 3480 m/s \\
    $c_{12}$ & 3.2 GPa & $T_3$ & 3.6 & $c_{s3}$ & 3520 m/s \\
    $c_{13}$ & 1.2 GPa & $K_f$ & 2.5 GPa & $c_{ps1}$ & 1030 m/s \\
    $c_{33}$ & 53.4 GPa & $\rho_f$ & 1040 kg/m$^3$ & $c_{ps3}$ & 746 m/s \\
    $c_{55}$ & 26.1 GPa & $\eta$ & $10^{-3}$ kg/m$\cdot$s & $\tau_{d1}$ & 5.95 $\mu$s \\
    $\phi$ & 0.2 & & & $\tau_{d3}$ & 1.82 $\mu$s \\
    \bottomrule
  \end{tabular}
\end{center}
\end{table}

All test cases are run with the orthotropic, transversely isotropic
sandstone of Table \ref{tab:matprops}, at a frequency of 10~kHz.  The
computational domain for each case is a cube with its center at the
origin, discretized with an equal number of cells in each direction;
for most cases the sides of the cube are aligned with the global
computational axes, but the grid is rotated for some cases to provide
a simple test of the mapped grid capabilities of the code.  For the fast P
wave and both S waves, the edge length of the domain is one wavelength
of the solution (computed as $2\pi/|\real{k}|$ for complex $k$), and
the total simulation time is 1.25 periods of the plane wave.  For the
slow P wave, the edge length is one decay length of the wave, computed
as $1/|\imag{k}|$, which is substantially less than one wavelength
even at this high frequency, and the total simulation time is set to
1.25 times the time for a fast P wave in the material $1$-direction to
cross the domain.  The simulation time step is
chosen so that the global maximum CFL number is 0.9.  For all cases,
boundary conditions are implemented by filling the ghost cells with
the true solution evaluated at cell centroids.  Limiting is not used
for any of the tests in this section in order to avoid obscuring the
convergence behavior of the underlying wave propagation algorithm.

\begin{table}
  \caption[Three-dimensional test cases]{Summary of plane wave test
    cases in three dimensions.  Within each group, cases are ordered
    by wave speed, fastest first.  The
    components of the $\Bell$ vector are given in grid axes.}
  \label{tab:plane-wave-cases-3d}
  \begin{center}
    \begin{tabular}{cccccccc}
      \toprule
      & \multicolumn{3}{c}{Grid axes} & \multicolumn{3}{c}{Material axes} & \\
      \cmidrule(r){2-4} \cmidrule(r){5-7}
      Cases & Yaw & Pitch & Roll & Yaw & Pitch & Roll & $\Bell$ vector \\
      \midrule
      0-3 & $0^\circ$ & $0^\circ$ & $0^\circ$ & $0^\circ$ & $0^\circ$ & $0^\circ$ & $(1,0,0)$ \\
      4-7 & $0^\circ$ & $0^\circ$ & $0^\circ$ & $0^\circ$ & $0^\circ$ & $0^\circ$ & $(0,0,1)$ \\
      8-11 & $30^\circ$ & $20^\circ$ & $10^\circ$ & $0^\circ$ & $0^\circ$ & $0^\circ$ & $(1,0,0)$ \\
      12-15 & $30^\circ$ & $20^\circ$ & $10^\circ$ & $0^\circ$ & $0^\circ$ & $0^\circ$ & $(0,1,0)$ \\
      16-19 & $30^\circ$ & $20^\circ$ & $10^\circ$ & $0^\circ$ & $0^\circ$ & $0^\circ$ & $(0,0,1)$ \\
      20-23 & $0^\circ$ & $0^\circ$ & $0^\circ$ & $30^\circ$ & $20^\circ$ & $10^\circ$ & $(1,0,0)$ \\
      24-27 & $0^\circ$ & $0^\circ$ & $0^\circ$ & $30^\circ$ & $20^\circ$ & $10^\circ$ & $(0,1,0)$ \\
      28-31 & $0^\circ$ & $0^\circ$ & $0^\circ$ & $30^\circ$ & $20^\circ$ & $10^\circ$ & $(0,0,1)$ \\
      32-35 & $0^\circ$ & $0^\circ$ & $0^\circ$ & $0^\circ$ & $0^\circ$ & $0^\circ$ &
      $\left(\frac{1}{\sqrt{3}}, \frac{1}{\sqrt{3}}, \frac{1}{\sqrt{3}} \right)$ \\
      \bottomrule
    \end{tabular}
  \end{center}
\end{table}

Table \ref{tab:plane-wave-cases-3d} lists the plane wave cases by
groups of four.  Within each group, each wave family is tested, in
decreasing order of speed --- the first case of each group is the fast
P wave, and the last is the slow P wave.  Cases 5 and 6 are shear
waves in the material's plane of isotropy, so their polarizations must
be explicitly specified; case 5 is polarized with its solid velocity
in the $x$ direction, while case 6 is polarized in the $y$ direction.
Cases 0-7 are the
simplest, with neither the grid nor the material principal axes
rotated from the global computational axes.  Note that these cases
only propagate waves in the $x$ and $z$ directions; the $y$ direction
would be redundant because the $x$-$y$ plane is the material's 1-2
plane, in which it is isotropic.  Cases 8-19 provide a basic test of
the mapped grid capabilities of the simulation code --- specifically
handling of grid interfaces that are not parallel to the
global coordinate planes --- while cases 20-31 test correct handling
of principal material directions that are not aligned with the global
axes.  The rotation matrix transforming from the grid or material axes
in cases 8-31 to the global $xyz$ axes is $\balpha = \B{R}_z(\psi)
\B{R}_y(-\phi) \B{R}_x(\theta)$, where $\psi$ is the yaw angle listed
in the table, $\phi$ is pitch, $\theta$ is roll, and
$\B{R}_\xi(\delta)$ is the elementary rotation matrix that rotates
counterclockwise by an angle $\delta$ about the $\xi$ axis.  All of
the above cases examine only waves propagating in the direction of one
of the grid axes, giving variation only in one grid direction for
dimensional splitting.  Cases 32-35, however, send waves
propagating obliquely through the grid, in order to see the full
effect of dimensional splitting on accuracy in three dimensions.

Tables \ref{tab:plane-wave-results-3d-gridaligned} and
\ref{tab:plane-wave-results-3d-nonaligned} list the results of
these convergence studies.  The 1-norm and max-norm errors in these
tables are normalized by the corresponding grid norm of the true solution.  The convergence behavior of the
three-dimensional code is exactly what would be expected from the
two-dimensional viscous high-frequency results of~\refcartpaper{} and
the formal order of accuracy of dimensional splitting: all
cases display second-order convergence in both the 1-norm and the
max-norm, except for cases 32-35, which involve waves not propagating
straight along the grid directions, and display first-order
convergence in both norms.  Because the solution is always
well-resolved on the fine grid, the error values are quite small in all cases; on the
$200^3$ grid, they only reach as high as 1.5\% relative error in the
max-norm for the off-axis cases, and 0.058\% relative error in the
max-norm for the grid-aligned cases.

\begin{table}
  \caption{Convergence results for the grid-aligned cases (numbers
    0-31) of Table \ref{tab:plane-wave-cases-3d}.}
  \label{tab:plane-wave-results-3d-gridaligned}
  \begin{center}
    \begin{tabular}{cccccccc}
      \toprule
      & & \multicolumn{3}{c}{Convergence rate} & \multicolumn{2}{c}{Error on $200^3$ grid} \\
      \cmidrule(r){3-5} \cmidrule(r){6-7}
      & Error norm & Best & Worst & Mean & Best & Worst \\
      \midrule
      \multirow{2}{*}{Fast P} & 1-norm & 2.05 & 2.03 & 2.05 & $8.81 \times 10^{-5}$ & $1.73 \times 10^{-4}$ \\
      & Max-norm & 2.01 & 1.96 & 2.00 & $2.09 \times 10^{-4}$ & $4.56 \times 10^{-4}$ \\
      \midrule
      \multirow{2}{*}{Shear 1} & 1-norm & 2.04 & 2.03 & 2.04 & $2.15 \times 10^{-4}$ & $2.22 \times 10^{-4}$ \\
      & Max-norm & 2.05 & 1.94 & 2.03 & $4.67 \times 10^{-4}$ & $5.86 \times 10^{-4}$ \\
      \midrule
      \multirow{2}{*}{Shear 2} & 1-norm & 2.04 & 2.03 & 2.04 & $2.19 \times 10^{-4}$ & $2.46 \times 10^{-4}$ \\
      & Max-norm & 2.05 & 1.94 & 2.01 & $4.73 \times 10^{-4}$ & $5.80 \times 10^{-4}$ \\
      \midrule
      \multirow{2}{*}{Slow P} & 1-norm & 2.02 & 2.02 & 2.02 & $9.50 \times 10^{-7}$ & $3.46 \times 10^{-6}$ \\
      & Max-norm & 1.95 & 1.83 & 1.92 & $1.57 \times 10^{-6}$ & $5.18 \times 10^{-6}$ \\
      \bottomrule
    \end{tabular}
  \end{center}
\end{table}

\begin{table}
  \caption{Convergence results for the non-grid-aligned cases (numbers
    32-35) of Table \ref{tab:plane-wave-cases-3d}.}
  \label{tab:plane-wave-results-3d-nonaligned}
  \begin{center}
    \begin{tabular}{cccc}
      \toprule
      Wave family & Error norm & Convergence rate & Error on $200^3$ grid \\
      \midrule
      \multirow{2}{*}{Fast P} & 1-norm & 1.01 & $4.19 \times 10^{-3}$ \\
      & Max-norm & 1.01 & $5.67 \times 10^{-3}$ \\
      \midrule
      \multirow{2}{*}{Shear 1} & 1-norm & 1.01 & $5.95 \times 10^{-3}$ \\
      & Max-norm & 0.91 & $1.26 \times 10^{-2}$ \\
      \midrule
      \multirow{2}{*}{Shear 2} & 1-norm & 1.01 & $6.99 \times 10^{-3}$ \\
      & Max-norm & 0.93 & $1.49 \times 10^{-2}$ \\
      \midrule
      \multirow{2}{*}{Slow P} & 1-norm & 1.00 & $1.26 \times 10^{-4}$ \\
      & Max-norm & 0.91 & $5.87 \times 10^{-4}$ \\
      \bottomrule
    \end{tabular}
  \end{center}
\end{table}

\subsection{Test of the revised limiter}
\label{sec:elimiter-test}

Having verified the numerical solution process without a limiter
present, it is now time to verify that the 
$\bE$-limiter defined in Section \ref{sec:elimiter} in fact
correctly limits shear waves on
non-rectilinear mapped grids.  To do this, a specially constructed
grid mapping is used, in order to make the two polarizations of
shear wave produced by the Riemann solver switch places when sorted in
order of speed.  While this example is somewhat contrived,
such a situation could easily happen by accident when using a more
realistic mapped grid in combination with an orthotropic material.

For this problem, the computational domain is the cube
$[-1,1]^3$, and the mapping from computational coordinates $\xi_1$,
$\xi_2$, $\xi_3$ to physical coordinates $x$, $y$, $z$ is
\begin{equation} \label{eq:tilt-map}
  x = \xi_1 L/2, \quad 
  y = \xi_2 L/2, \quad
  z = \begin{cases}
    (\xi_3 + \sigma \xi_1 \xi_3^3) L/2, \quad & \xi_3 < 0 \\
    (\xi_3 + \sigma \xi_2 \xi_3^3) L/2, \quad & \xi_3 \ge 0.
  \end{cases}
\end{equation}
This is a $\mathcal{C}^2$ map whose $\xi_3$ grid surfaces are tilted
in the $x$ direction for $\xi_3 < 0$, but in the $y$ direction for
$\xi_3 > 0$.  The slope parameter $\sigma$ is set to 0.1.
Using the orthotropic sandstone of Table \ref{tab:matprops}, with the
material 1-2-3 axes coinciding with the $x$-$y$-$z$ axes, this means
that the $x$-direction shear waves found by the Riemann solver will be
the faster polarization for $z < 0$, but the $y$-direction waves will
be faster for $z > 0$.  If the common convention of sorting the waves
by their speeds is used, then as discussed in Section
\ref{sec:elimiter}, the conventional wave strength ratio
\eqref{eq:3dlim-theta-classical} can be expected to cause
difficulties.

\begin{table}
  \caption{Comparison of performance of classical limiter wave strength
    ratio and $\bE$-limiter for plane wave case 5 of Table
    \ref{tab:plane-wave-cases-3d}, using the tilted grid map of
    \eqref{eq:tilt-map}.}
  \label{tab:elimiter-comparison}
  \begin{center}
    \begin{tabular}{ccccc}
      \toprule
      \multirow{2}{*}[-3.5pt]{\parbox{0.55in}{\centering Strength ratio}}
      & & \multicolumn{3}{c}{Relative error on grid of size} \\
      \cmidrule{3-5}
      & Error norm & $50^3$ & $100^3$ & $200^3$ \\
      \midrule
      \multirow{2}{*}{Classical} & 1-norm & $5.97 \times 10^{-3}$ & $1.85 \times 10^{-3}$ & $6.43 \times 10^{-4}$ \\
      & Max-norm & $4.21 \times 10^{-2}$ & $2.20 \times 10^{-2}$ & $1.16 \times 10^{-2}$ \\
      \midrule
      \multirow{2}{*}{$\bE$-limiter} & 1-norm & $3.97 \times 10^{-3}$ & $1.29 \times 10^{-3}$ & $5.05 \times 10^{-4}$ \\
      & Max-norm & $1.58 \times 10^{-2}$ & $5.77 \times 10^{-3}$ & $2.14 \times 10^{-3}$ \\
      \bottomrule
    \end{tabular}
  \end{center}
\end{table}

To show the difference between the two wave strength ratios to
greatest effect, this grid is used to simulate test case 5 of Table
\ref{tab:plane-wave-cases-3d} --- a shear wave propagating in the $+z$
direction, and polarized in the $x$ direction.  The length parameter
$L$ in the mapping is set to one wavelength of the wave; all other
parameters are identical to the rectilinear grid case.  Table
\ref{tab:elimiter-comparison} compares the results using the new wave
strength ratio of \eqref{eq:3dlim-theta} and the old ratio
\eqref{eq:3dlim-theta-classical}, both using the MC limiter function.
The $\bE$-limiter reduces the error in both the 1-norm and the
max-norm; while the 1-norm error reduction is modest, the reduction of
the max-norm error is quite substantial, up to a factor of five on the
finest grid.

\subsection{Demonstration problem}
\label{sec:3d-demo}

\begin{figure}
  \begin{center}
    \includegraphics[width=0.5\textwidth]{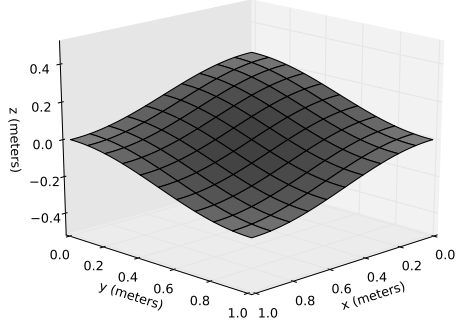}
    \caption[Brine-sandstone interface surface for demonstration
      problem]{Brine-sandstone interface surface for demonstration
      problem. \label{fig:undulate-surfplot}}
  \end{center}
\end{figure}

To demonstrate the numerical methods developed 
here, and the capabilities of the code implementing them, this
section discusses the simulation of an acoustic pulse in brine
striking an undulating bed of orthotropic layered sandstone.  The
surface of the bed is defined by
\begin{equation}
  z_{\text{int}}(x,y) = z_0 + H_x \cos \left( \frac{2 \pi x}{L_x}
  \right)
  + H_y \cos \left( \frac{2 \pi y}{L_y} \right),
\end{equation}
with the parameters $z_0$, $H_x$, $L_x$, $H_y$, and $L_y$ given in
Table \ref{tab:undulate-params}.  Figure \ref{fig:undulate-surfplot}
shows the surface.  Below this $z$ coordinate, the domain is composed
of the orthotropic sandstone of Table \ref{tab:matprops}; above, it is
composed of the brine from this sandstone.  The curved interface
between the two media is incorporated into the model using a
mapped grid.  At every point within the sandstone, the material's
plane of isotropy (the plane of the principal 1-2 axes) is parallel to
the tangent plane of the surface above, in order to simulate a bed
that has been folded, or deposited on a pre-existing uneven
surface.  In the simulation code, these variable principal axes
are implemented by assigning constant material principal directions
to each cell, equal to the directions evaluated at the cell
centroid.  The interface is taken to have open pores ($\eta_d = 1$ in
interface condition \eqref{eq:poro-fluid-physical}), and the incoming
acoustic pulse propagates straight downward in the $-z$ direction.
This problem exercises almost all of the capabilities of the
three-dimensional code --- it involves mapped grids, an orthotropic
material with a variable principal direction, and a fluid-poroelastic
interface.

\begin{table}
  \caption[Parameters for demonstration problem]{Surface and mapping
    parameters for demonstration problem.  Note that the domain has
    horizontal dimensions $L_x/2$ by $L_y/2$.}
  \label{tab:undulate-params}
  \begin{center}
    \begin{tabular}{lrrrrrrrrrrrr}
      \toprule
      Surface parameter & $z_0$ & $L_x$ & $L_y$ & $H_x$ & $H_y$ \\
      \cmidrule(r){1-6}
      Value & 0~m & 2~m & 2~m & $\frac{3 L_x}{16 \pi}$ & $\frac{3
        L_y}{16 \pi}$ \\
      \midrule
      Mapping parameter & $z_{\text{bot}}$ & $z_{\text{top}}$ &
      $\xi_{\text{bot}}$ & $\xi_{\text{int}}$ & $\xi_{\text{top}}$ &
      $r_{\text{bot}}$ & $r_{\text{top}}$ \\
      \cmidrule(r){1-8}
      Value & $-1$~m & 0.5~m & 0.15 & 0.6 & 0.9 &
      $\frac{2}{\xi_{\text{bot}}}$ & $\frac{2}{1-\xi_{\text{top}}}$ \\
      \bottomrule
    \end{tabular}
  \end{center}
\end{table}

The grid mapping for this problem is defined so that one of
computational coordinate surfaces follows the interface, with the rest
of the map chosen as a compromise between simplicity, smoothness, and
the ability to have a flat grid plane at a useful distance below the
interface to output slices of the solution for later plotting.
The computational domain is the unit cube $[0,1]^3$; in the $xy$
plane, the problem's symmetry allows the physical domain to be chosen
as one quarter of a periodic tile of the surface, $[0, L_x/2] \times
[0, L_y/2]$.  The grid mapping function in the horizontal axes is a
simple scaling, $x := \xi_1 L_x/2$ and $y := \xi_2 L_y/2$, while the
mapping function for the $z$ coordinate is defined in terms of the
horizontal physical coordinates $x$ and $y$, and the vertical
computational coordinate $\xi_3$, as
\begin{equation} \label{eq:z-undulate}
  z := \begin{cases}
    z_{\text{bot}} + \frac{z'_{\text{bot}}}{r_{\text{bot}}} \sinh
    r_{\text{bot}} (\xi_3 - \xi_{\text{bot}}),  & \xi_3 <
    \xi_{\text{bot}} \\
    z_{\text{bot}} + z'_{\text{bot}} (\xi_3 - \xi_{\text{bot}})
    + a_{\text{bot}}(x,y) \, b(\xi_3; \xi_{\text{bot}}),  & \xi_{\text{bot}} \le \xi_3 <
    \xi_{\text{int}} \\
    z_{\text{top}} + z'_{\text{top}} (\xi_3 - \xi_{\text{top}})
    + a_{\text{top}}(x,y) \, b(\xi_3; \xi_{\text{top}}),  & \xi_{\text{int}} \le \xi_3 <
    \xi_{\text{top}} \\
    z_{\text{top}} + \frac{z'_{\text{top}}}{r_{\text{top}}} \sinh
    r_{\text{top}} (\xi_3 - \xi_{\text{top}}),  & \xi_3 \ge
    \xi_{\text{top}},
  \end{cases}
\end{equation}
where the derived quantities $z'_{\text{bot}}$, $z'_{\text{top}}$,
$a_{\text{bot}}(x,y)$, $a_{\text{top}}(x,y)$, and $b(\xi_3; \xi^*)$ in
the mapping function are defined by
\begin{equation}
\begin{gathered}
  \begin{aligned}
  z'_{\text{bot}} &:= \frac{z_0 - H_x - H_y -
    z_{\text{bot}}}{\xi_{\text{int}} - \xi_{\text{bot}}}
  &z'_{\text{top}} &:= \frac{z_{\text{top}} - z_0 - H_x -
    H_y}{\xi_{\text{top}} - \xi_{\text{int}}} \\
  a_{\text{bot}}(x,y) &:= z_{\text{int}}(x,y) - (z_0 - H_x - H_y)
  &a_{\text{top}}(x,y) &:= z_{\text{int}}(x,y) - (z_0 + H_x + H_y)
  \end{aligned}\\
  b(\xi_3; \xi^*) := \frac{1}{2} \left(\sqrt{1 + 8 \left( \frac{\xi_3 -
        \xi^*}{\xi_{\text{int}} - \xi^*} \right)^2} - 1 \right).
\end{gathered}
\end{equation}
The values of the mapping parameters are given in Table
\ref{tab:undulate-params}.

\begin{figure}
  \begin{center}
    \subfloat[Back $yz$ face ($x = 0$)]
             {\includegraphics[width=0.3\textwidth]{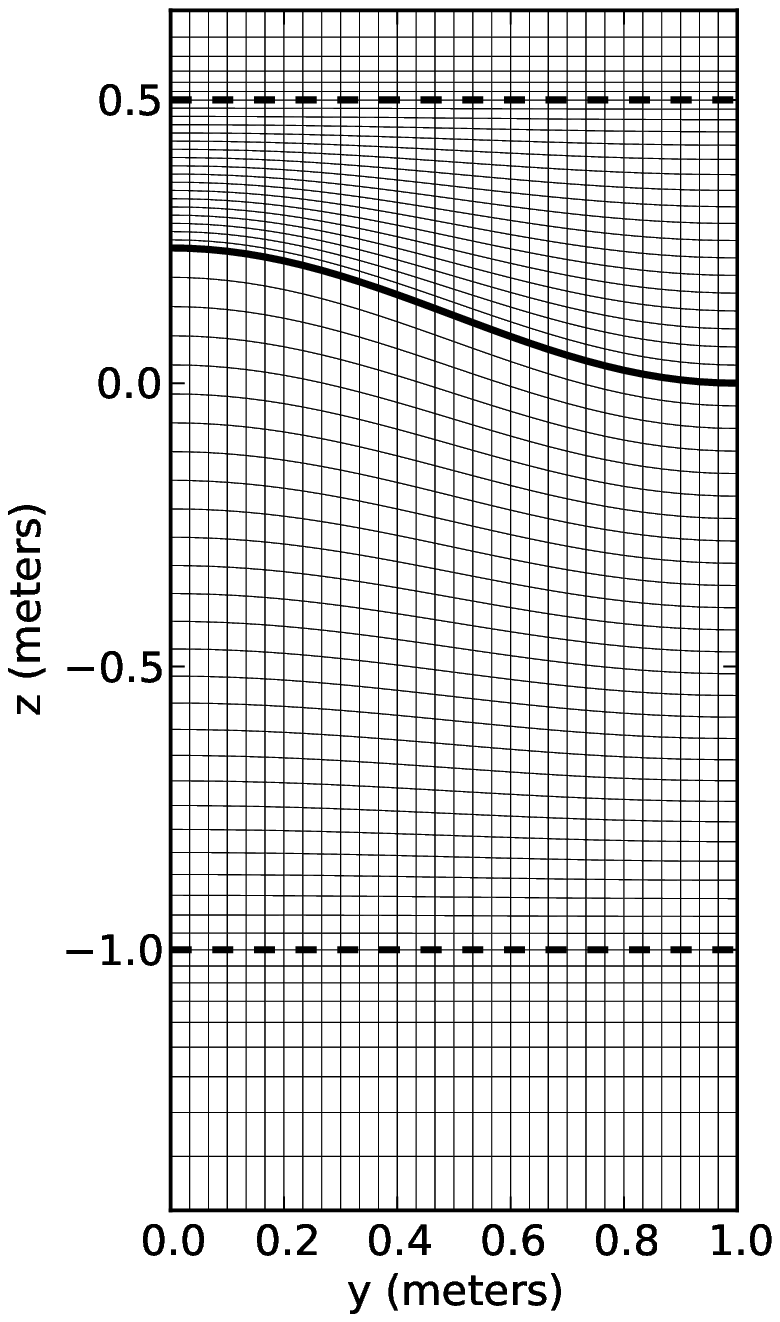}}
    \hspace{0.5in}
    \subfloat[Front $yz$ face ($x = L_x/2$)]
             {\includegraphics[width=0.3\textwidth]{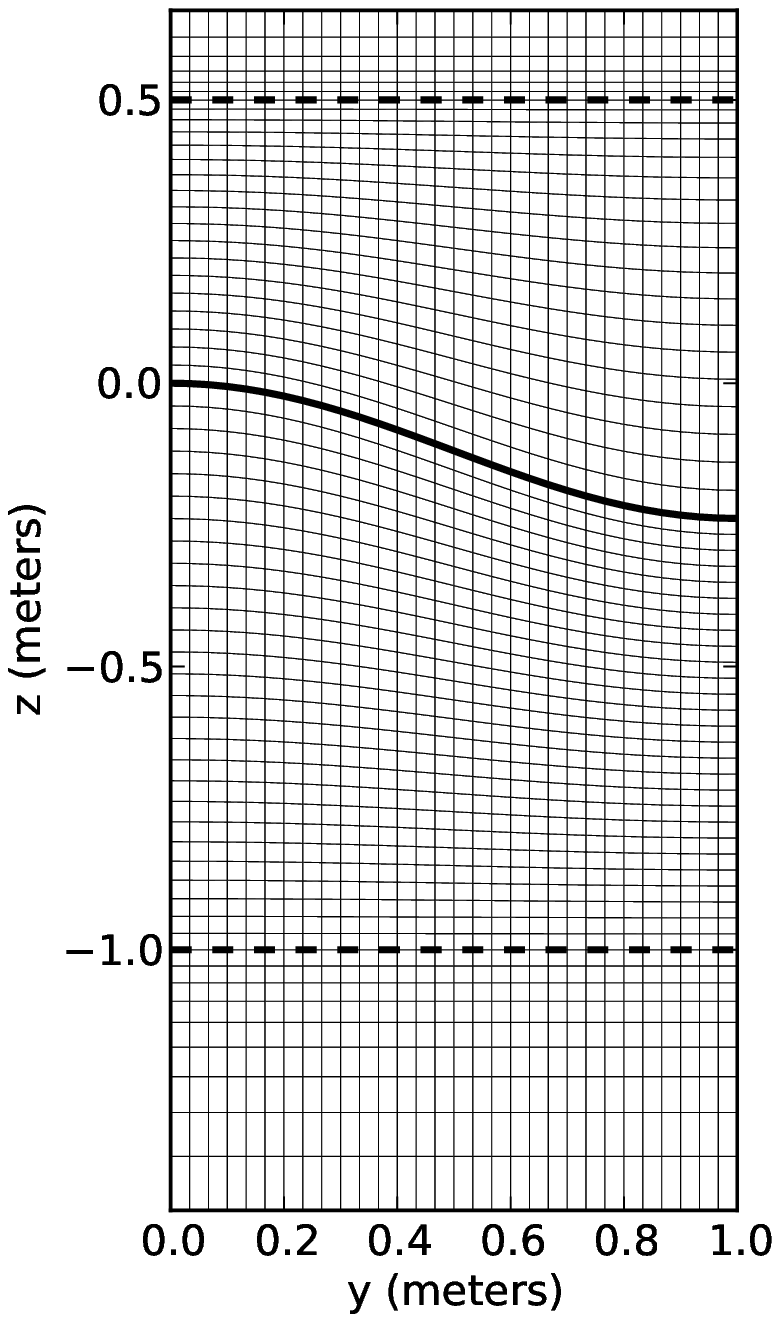}}
    \caption[Grid mapping for demonstration problem]{Side views of the
      mapped grid used for the demonstration problem, displayed on a
      $30 \times 30 \times 60$ grid.  The heavy solid line indicates
      the sandstone-brine interface; the heavy dashed lines mark the
      boundaries of the regions where the grid is stretched to move
      the boundaries outward.  Because $H_y = H_x$ and $L_y = L_x$,
      plots of the $xz$ faces would look
      identical. \label{fig:undulate-map}}
  \end{center}
\end{figure}

The intent of the mapping \eqref{eq:z-undulate} is to provide uniform
grid spacings in the shortest grid columns above and below the
interface, in order to prevent cells in these columns from being any
smaller than necessary.  The parameters $z'_{\text{top}}$ and
$z'_{\text{bot}}$ are rates of change of $z$ with respect to $\xi_3$
in these shortest columns.  Coordinates $z_{\text{top}}$,
$\xi_{\text{top}}$, $z_{\text{bot}}$, and $\xi_{\text{bot}}$ are the
boundaries in physical and computational space of the part of the
domain where the solution is considered ``interesting.''  Beyond these
coordinates, the hyperbolic sine term smoothly stretches the grid in
the vertical direction to move the boundaries of the computational
domain further away from the interface, with the intent of improving
the performance of the non-reflecting boundary conditions.  The
stretching will tend to blur the solution in the elongated cells, but
this is acceptable because a high-quality solution is not required
more than a little distance outside the region between
$z_{\text{bot}}$ and $z_{\text{top}}$.  Note that for $\xi_3 \le
\xi_{\text{bot}}$ and $\xi_3 \ge \xi_{\text{top}}$, surfaces of
constant $\xi_3$ are horizontal planes.  This is convenient for
outputting a horizontal slice of the solution just below
$z_{\text{bot}}$ for plotting.  In the middle portion of the grid, the
mapping is designed to have a continuous first derivative everywhere
except at the sandstone-brine interface, in order to avoid any
possible spurious internal reflections that might be caused by a
nonsmooth mapping, and to improve accuracy in general.  The mapping is
allowed to be nonsmooth at the sandstone-brine interface because the
second-order correction term is omitted there in any case, as
in~\refmappaper{} --- accuracy
will already be degraded there, and requiring the mapping to be smooth
at the interface would result in greater cell size variation
elsewhere.  Figure \ref{fig:undulate-map} shows the resulting grid.
Note that the closely-spaced cells above the interface are not
problematic for stability --- since they are in the brine, not the
sandstone, the wave speed within them is the acoustic wave speed of
1550~m/s, whereas the fast P wave speed in the sandstone is always at
least 5260~m/s.  The smallest vertical dimension of the cells above
the interface is still over half the smallest vertical dimension below
it, so stability is restricted by the fast P wave in the sandstone,
not the acoustic wave in the brine.

By symmetry, the boundary conditions at the lateral faces of the
domain are set to be reflective --- for the faces parallel to the $yz$
plane, ghost cells are set to the value of the adjacent cell in the
computational domain but with $\tau_{xz}$, $\tau_{xy}$, $v_x$, and
$q_x$ negated, while for the faces parallel to the $xz$ plane,
$\tau_{yz}$, $\tau_{xy}$, $v_y$, and $q_y$ are negated.
Non-reflecting boundary conditions at the top and bottom face are
implemented using zero-order extrapolation, with the elongation of the
grid mapping near the top and bottom boundaries used to move these
boundaries further away from the interface.  Moving the boundaries
further away allows the waves generated at the sandstone-brine
interface to have angles of incidence closer to normal, which reduces
reflections from this simple approach, and also postpones the arrival
of waves at the boundary.  The initial state is
set to zero everywhere, except for the incoming plane wave, which is
defined by its pressure field,
\begin{equation}
  p_{\text{in}}(x,y,z) := \begin{cases}
    0.5\,\text{Pa} \left( 1 + \cos \left(
  \frac{2\pi (z - z_{0,\text{wave}})}{\lambda_{\text{wave}}} \right)
  \right), \quad & |z - z_{0,\text{wave}}| < \lambda_{\text{wave}}/2\\
  0, \quad & \text{otherwise},
  \end{cases}
\end{equation}
where $z_{0,\text{wave}} = z_0 + H_x + H_y +
0.6\lambda_{\text{wave}}$, $\lambda_{\text{wave}} =
c_{\text{wave}}/f_{\text{wave}}$, $c_{\text{wave}}$ is the sound speed
in the brine (1550~m/s), and the fundamental frequency
$f_{\text{wave}}$ is 10~kHz.  To give a downward-propagating acoustic
wave, the vertical fluid velocity is set to $q_z =
-p_{\text{in}}/Z_f$, where $Z_f$ is the acoustic impedance of the
brine.  The total simulation time is 400~$\mu$s, and the dimensions of
the grid used are $300 \times 300 \times 600$ cells.  The MC limiter
is used with all waves, with the full energy inner product wave
strength ratio \eqref{eq:3dlim-theta}.

In order to obtain a solution more quickly, this demonstration problem
was run in parallel using \textsc{PetClaw}.  While there is a
significant memory overhead associated with the PETSc distributed
arrays used by \textsc{PetClaw}, as of this writing it is the only
\clawpack{} variant capable of running in parallel using dimensional
splitting, which made it the only practical option in terms of run
time for a large three-dimensional problem.  With \textsc{PetClaw} the
maximum memory footprint for the $300 \times 300 \times 600$ cell grid
was roughly 105~GB.  The problem was run on an
Amazon EC2 CR1 high-memory cluster compute node, with 16 MPI
processes; a total of 835 time steps were required, with a run time of
25.5 hours without viscosity included or 26.5 hours with viscosity,
giving an aggregate throughput of roughly 30,000 cell time steps per
CPU-second on the Intel Xeon E5-2670 CPUs used on this machine.

\begin{figure}
  \begin{center}
    \subfloat[Energy density ($\text{J}/\text{m}^3$) \label{fig:undulate-energy-invisc}]
             {\includegraphics[width=0.49\textwidth]{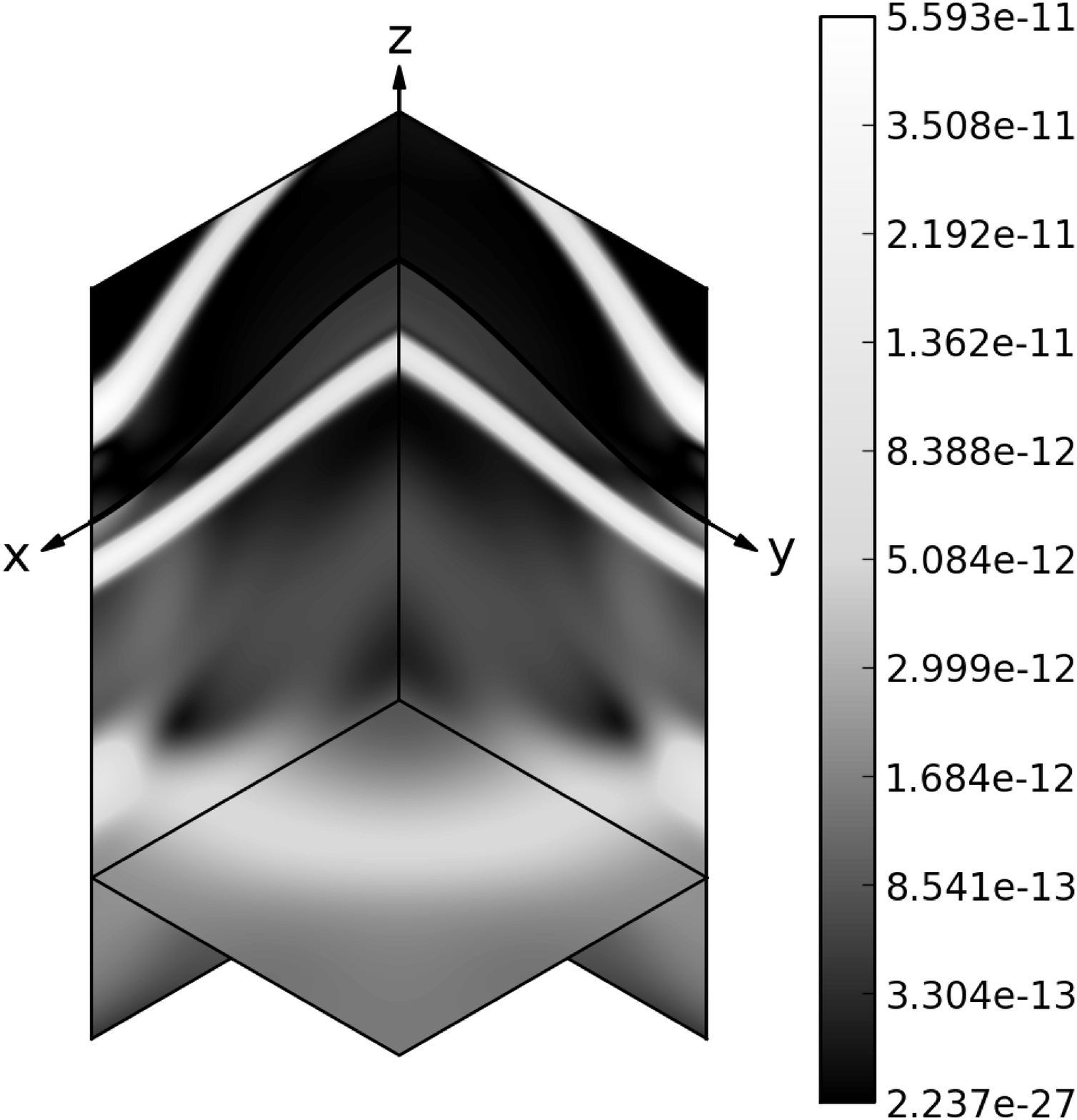}}
    \subfloat[Fluid pressure (Pa) \label{fig:undulate-p-invisc}]
             {\includegraphics[width=0.49\textwidth]{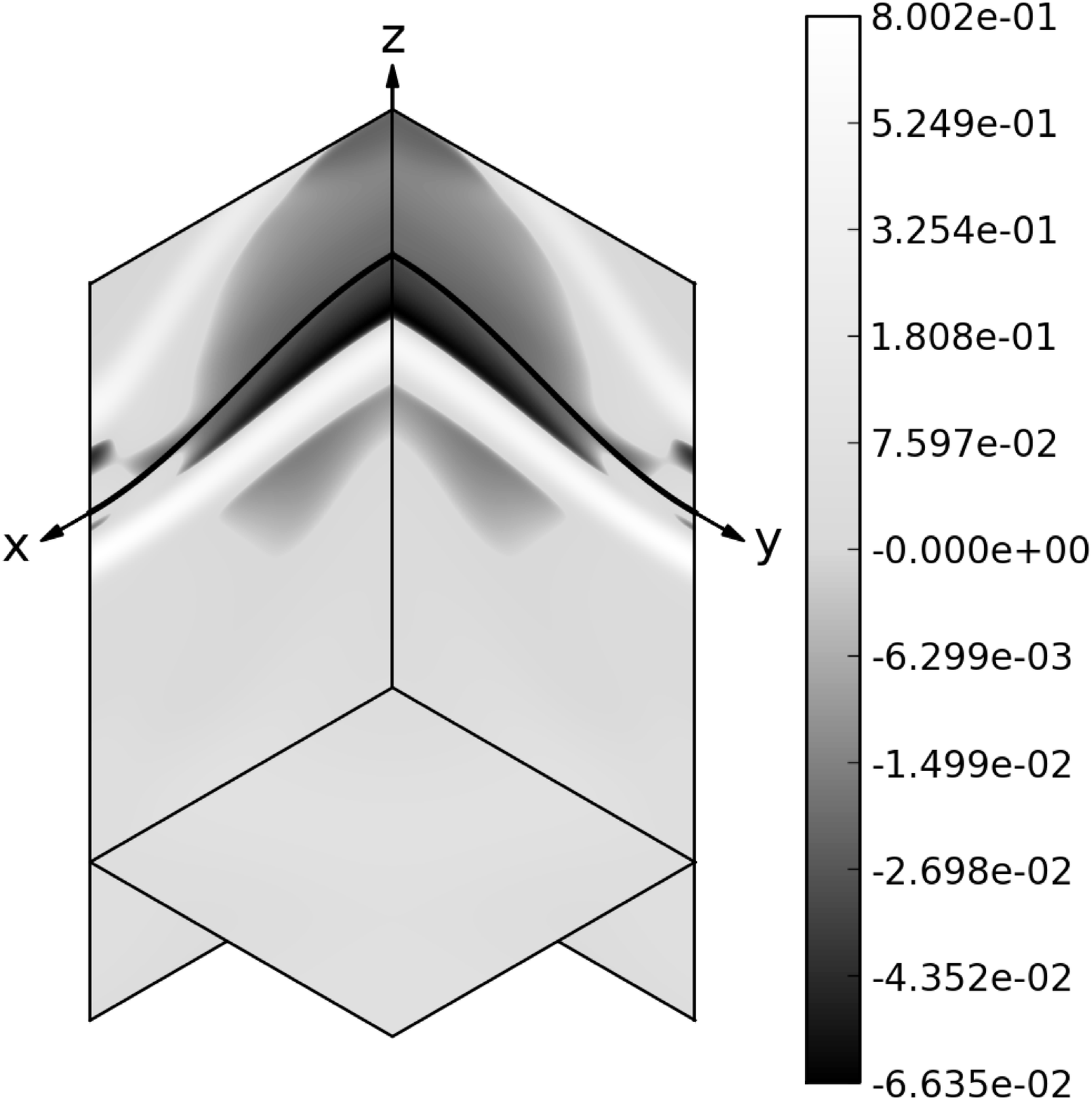}}\\
    \subfloat[Vertical direction normal stress $\tau_{zz}$ (Pa) \label{fig:undulate-tauzz-invisc}]
             {\includegraphics[width=0.49\textwidth]{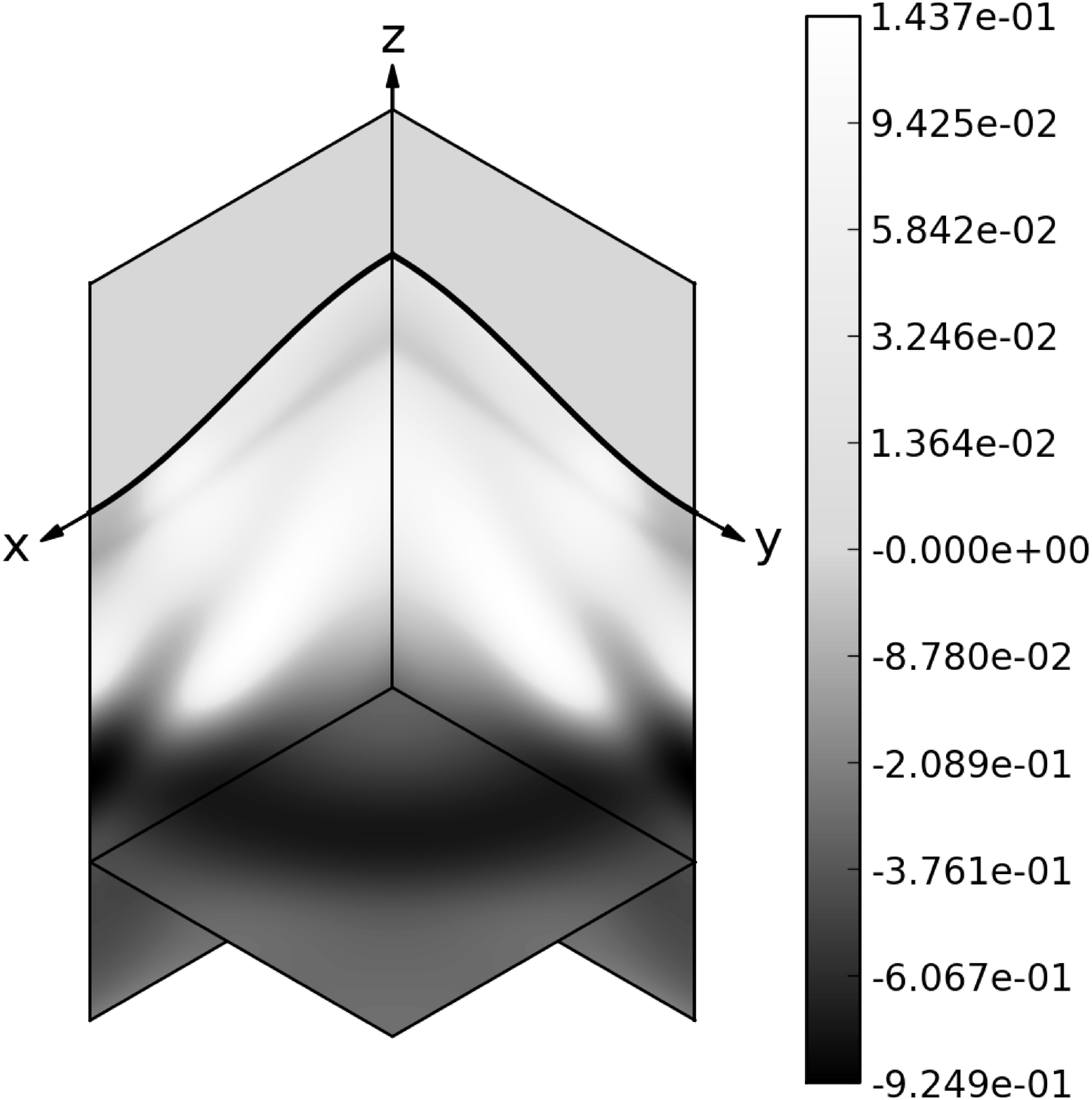}}
    \subfloat[Vertical direction solid velocity $v_z$ (m/s) \label{fig:undulate-vz-invisc}]
             {\includegraphics[width=0.49\textwidth]{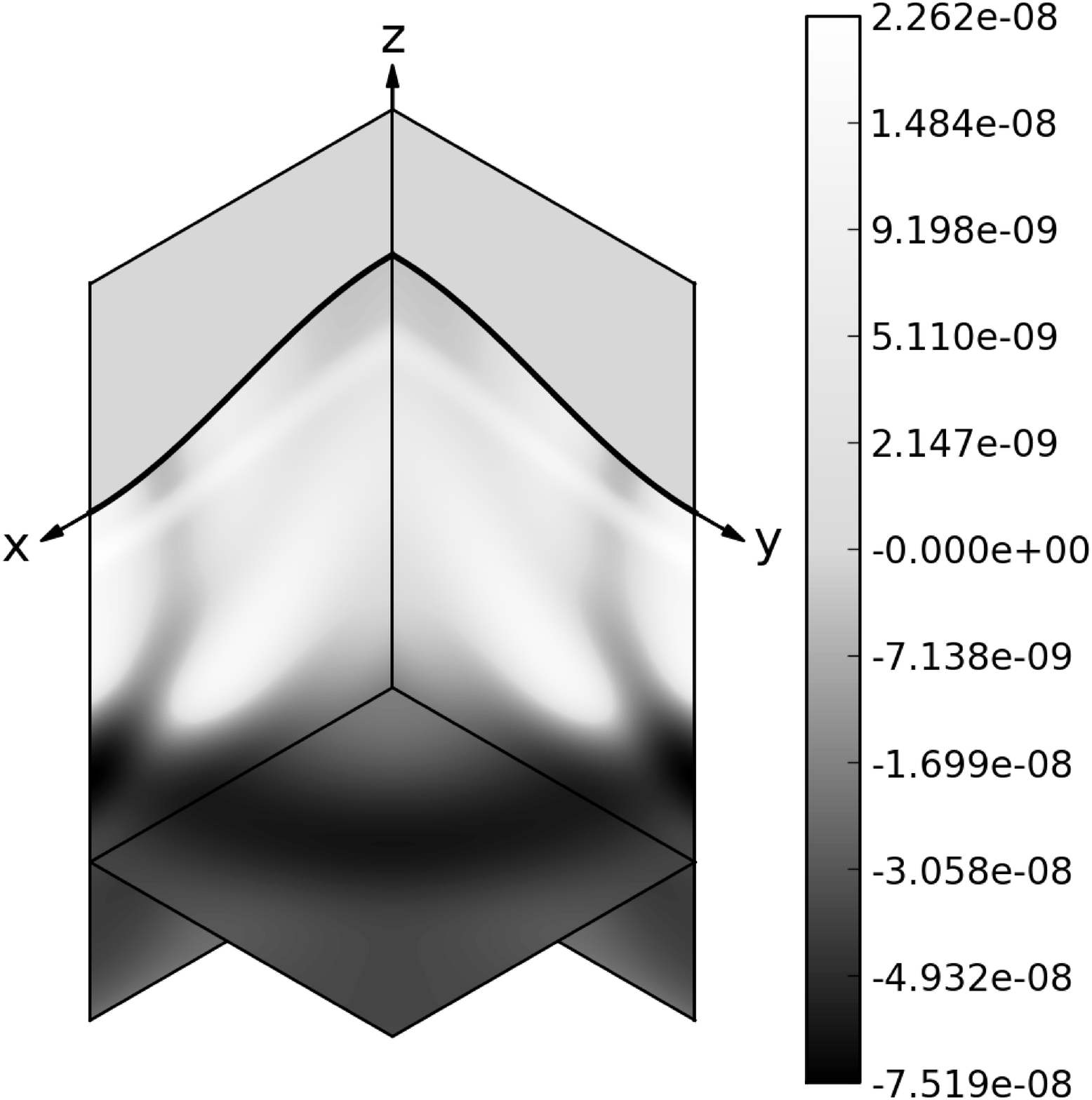}}
    \caption{Results for the demonstration problem at time
      399.9~$\mu$s, without viscosity. \label{fig:undulate-results-invisc}}
  \end{center}
\end{figure}

\begin{figure}
  \begin{center}
    \subfloat[Energy density ($\text{J}/\text{m}^3$) \label{fig:undulate-energy-visc}]
             {\includegraphics[width=0.49\textwidth]{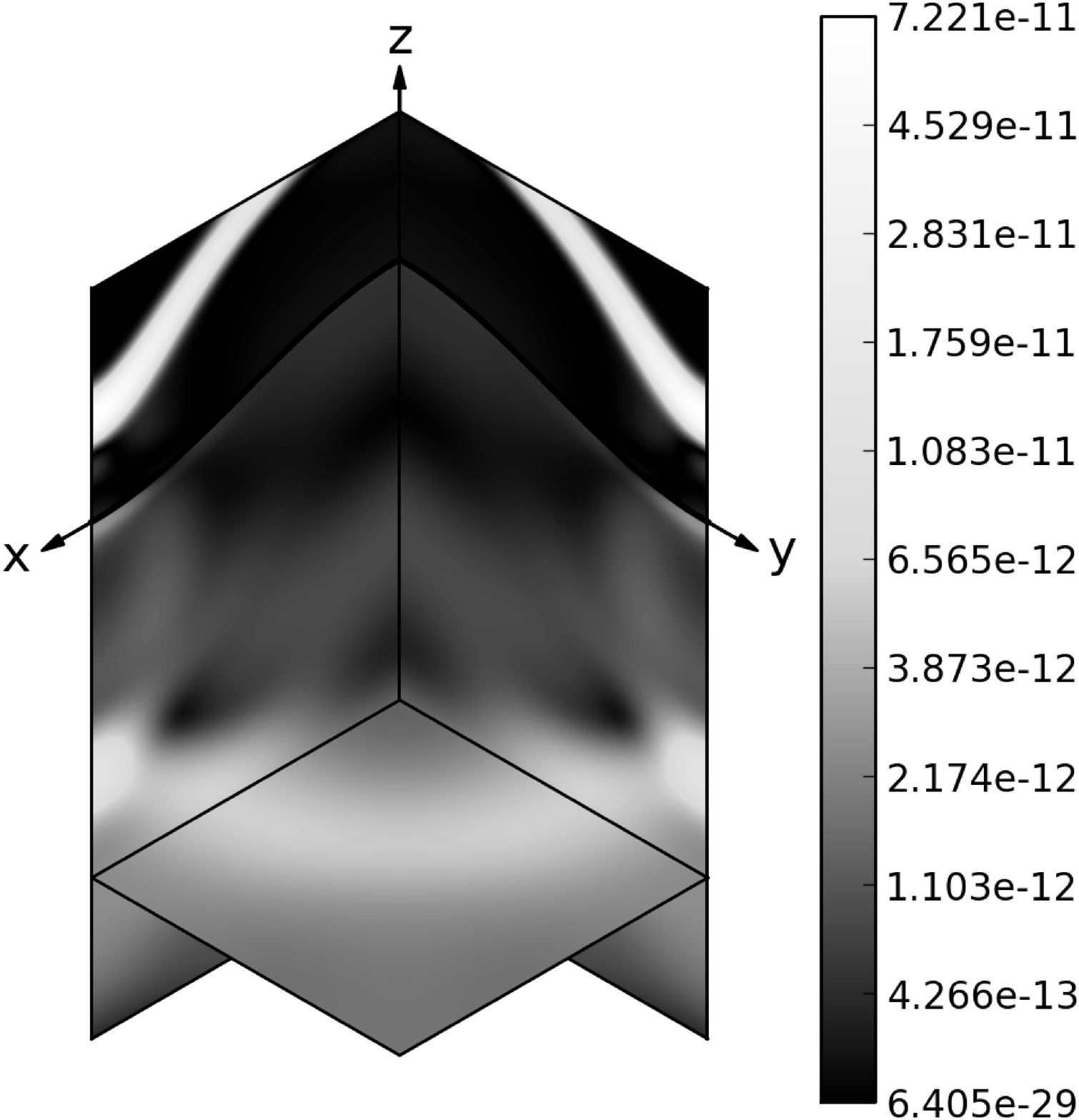}}
    \subfloat[Fluid pressure (Pa) \label{fig:undulate-p-visc}]
             {\includegraphics[width=0.49\textwidth]{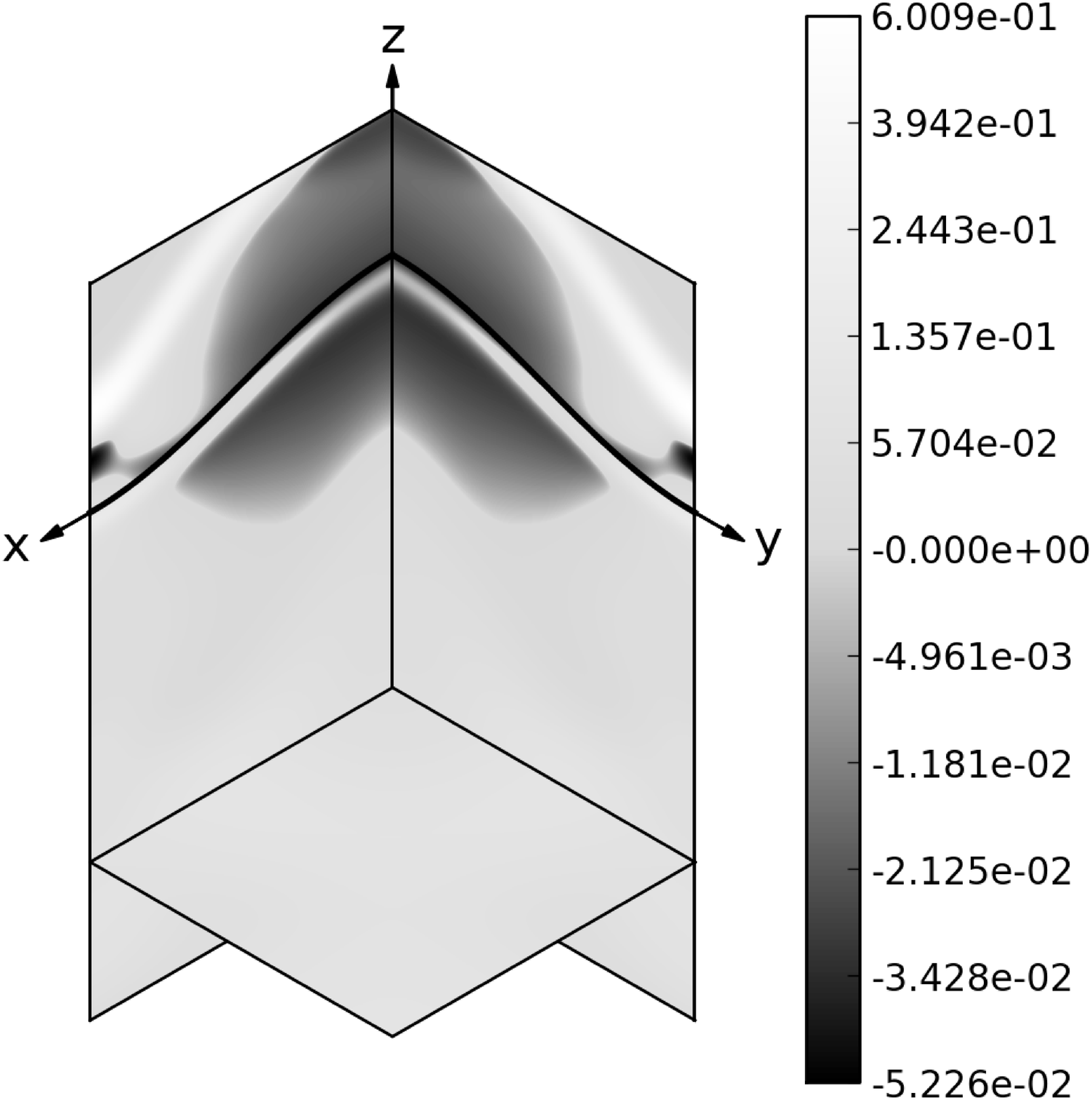}}\\
    \subfloat[Vertical direction normal stress $\tau_{zz}$ (Pa) \label{fig:undulate-tauzz-visc}]
             {\includegraphics[width=0.49\textwidth]{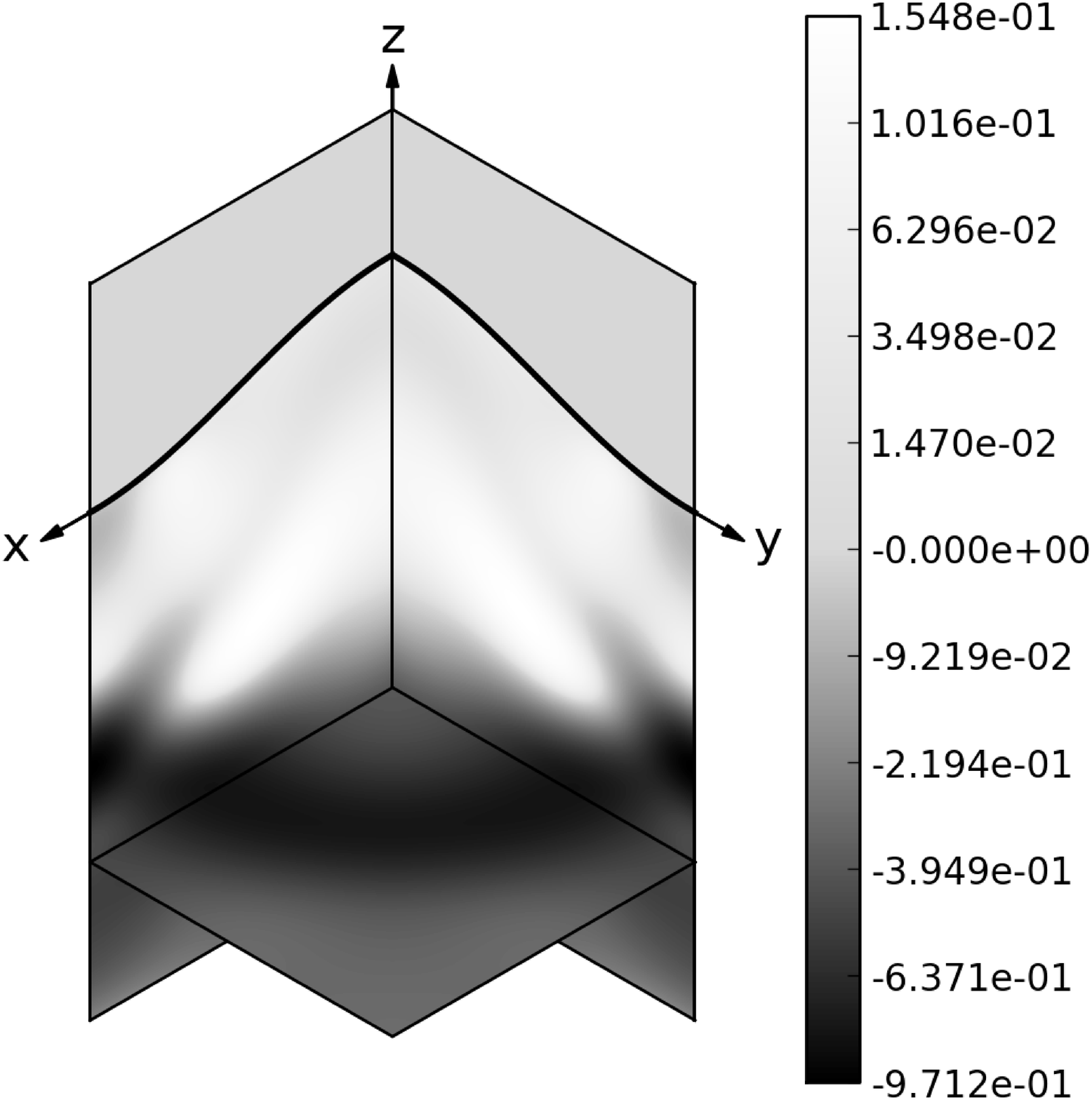}}
    \subfloat[Vertical direction solid velocity $v_z$ (m/s) \label{fig:undulate-vz-visc}]
             {\includegraphics[width=0.49\textwidth]{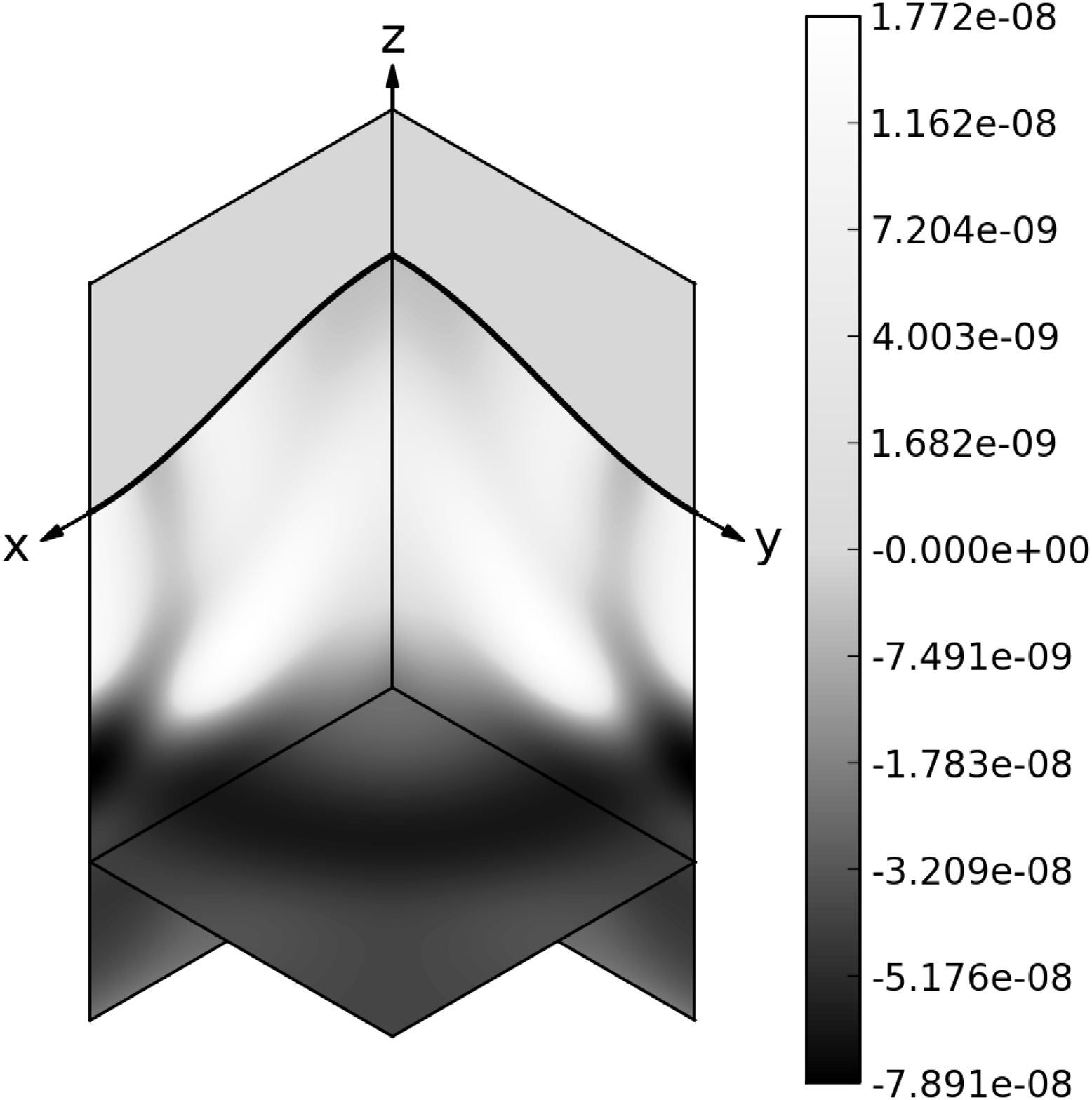}}
    \caption{Results for the demonstration problem at time
      399.9~$\mu$s, with viscosity. \label{fig:undulate-results-visc}}
  \end{center}
\end{figure}

Figures \ref{fig:undulate-results-invisc} and
\ref{fig:undulate-results-visc} show the solution at time
399.9~$\mu$s, the beginning of the final time step, with Figure
\ref{fig:undulate-results-invisc} showing the solution without
viscosity included and Figure \ref{fig:undulate-results-visc} showing
it with viscosity.  The plots show an isometric view of the
computational domain, with plots rendered on the $xz$ and $yz$ planes
and on the horizontal plane at the $z$ coordinate of the centroids of
the first layer of cells below $z = z_{\text{bot}}$; this is the
highest layer of cells in the sandstone whose centroids all lie in a
horizontal plane, so that a plot rendered on the surface defined by
these centroids is easy to interpret.  The plotted values on the $xz$
and $yz$ planes are generated by first projecting the values from the
centroids of the layer of cells on that side of the domain to the
appropriate plane, using the problem symmetry --- that is, for values
on the $xz$ plane, $\tau_{yz}$, $\tau_{xy}$, $v_y$, and $q_y$ are set
to zero, while for values on the $yz$ plane, $\tau_{xz}$, $\tau_{xy}$,
$v_x$, and $q_x$ are set to zero.  The locations associated with these
values for plotting purposes are the orthogonal projections of the
corresponding cell centroids onto the axis planes.  The problem
symmetry is also used to extend the computed solution to the lateral
corners of the computational domain -- values at $(x,y) = (0,0)$,
$(L_x/2, 0)$, $(0, L_y/2)$, and $(L_x/2, L_y/2)$ are obtained by copying the values at
the nearest cells, then setting all the shear stresses and
horizontal velocity components to zero.  In addition, since computing
the energy density requires knowing the material principal directions,
the principal directions on the $xz$ and $yz$ faces are computed at
the points on these planes associated with the projected values, not
the original cell centroids.  Values plotted in the horizontal plane
near the bottom of the domain are associated with the $z$ coordinate
of the cell centroids, which is $-1.0014$~m; this is the $z$
coordinate of the plane shown.

The solution of this demonstration problem is quite complex, but there
are a number of clearly recognizable features.  For the inviscid case,
Figure \ref{fig:undulate-energy-invisc}, which shows the energy
density, provides a broad view with most of the solution features
identifiable.  At the bottom, the light arc across the horizontal
slice and stretching up into the lower parts of the sides of the
domain is the initial fast P wave created when the acoustic wave
struck the peak of the sandstone.  The additional arc sweeping up and
inward from the intersection of the fast P wave with the domain edge
is the same initial fast P wave reflected off the boundary.  Further
inward toward the $z$ axis, the light diagonal bands are shear waves
originating from the acoustic wave striking the flanks of the
sandstone peak.  Just below the surface, the slow P wave is clearly
visible as a narrow bright band; because the pore structure is open at
the surface, a strong slow P wave is excited by the incident acoustic
wave.  Finally, above the surface, the incident wave has been
reflected and has already partially left the computational domain.
Figure \ref{fig:undulate-p-invisc} shows some numerical artifacts in
the low-pressure region at the very top of the domain, but these are
outside the designated area of interest for the problem.

For the viscous case, the solution is generally similar, but the slow
P wave is almost entirely suppressed by the viscous dissipation.  It
is, however, faintly visible as a band of increased pressure just
under the surface in Figure \ref{fig:undulate-p-visc}.  Comparing the
energy density plots of Figures \ref{fig:undulate-energy-invisc} and
\ref{fig:undulate-energy-visc}, the fast P and S waves are also
somewhat dissipated and slowed by viscosity, although the vertical
direction stress fields (Figures \ref{fig:undulate-tauzz-invisc} and
\ref{fig:undulate-tauzz-visc}) are hardly affected aside from the loss
of the slow P wave.

\section{Summary and future work}
\label{sec:conclusion}

This paper has covered the extension of the finite volume wave
propagation methods for poroelasticity developed in~\refmappaper{}
and~\refcartpaper{} to three dimensions.  Section \ref{sec:theory}
covered the development of a first-order linear hyperbolic system of
PDEs describing three-dimensional Biot theory at low frequencies.  An
energy density functional $\mathcal{E}$ was developed for the
three-dimensional system, and as in~\refcartpaper{} $\mathcal{E}$ was
found to be a strictly convex entropy function of the system in the
sense of Chen, Levermore, and
Liu~\cite{chen-levermore-liu:stiff-relaxation}.  Interface conditions
for coupling fluid and poroelastic media were also exhibited.

Section \ref{sec:mapped-fvm-3d} discussed the implementation of
high-resolution finite volume methods for fluid-poroelastic problems
on mapped grids in three dimensions.  The complications associated
with defining cell face normals on an arbitrary hexahedral grid were
discussed, and a technique was developed for defining suitable normal
vectors for a finite volume scheme.  The solution procedure for
poroelastic-fluid and poroelastic-poroelastic Riemann problems with
interface conditions developed in~\refmappaper{} was also extended to
three dimensions.  In addition, a new strength ratio for wave limiting
was developed for three-dimensional poroelasticity, which avoided the
problems with ambiguous shear wave polarization directions that would
otherwise be encountered; besides avoiding the inappropriate
suppression of higher-order terms that could be encountered with the
traditional wave strength ratio calculation, the new limiting approach
also gave a modest reduction in error for most cases when applied to
the cylindrical scatterer test problems of~\refmappaper{}.

With all the algorithmic pieces in place, Section \ref{sec:results}
applied the methods of Section \ref{sec:mapped-fvm-3d} to some test
problems to verify their effectiveness.  The first test problems were
simple plane waves, for which the numerical solution could easily be
compared to an analytical solution for the same problem.  Due to the
use of dimensional splitting, only first-order accuracy could be
achieved in the general case of waves propagating obliquely to the
grid, although when the wavevector was aligned with the grid axes
second-order convergence was achieved, consistent with previous
results~\cite{gil-ou:poro-2d-mapped, gil-ou-rjl:poro-2d-cartesian}.  A
special set of test problems was then run to demonstrate the
new $\bE$-limiter on a problem of the type it was developed for, where
the different polarizations of shear wave switch order in the Riemann
solution output; the $\bE$-limiter gave a modest reduction in 1-norm
error, and a substantial reduction --- up to a factor of five --- in
max-norm error.   A more complex demonstration problem involving an
acoustic wave in brine striking a periodically undulating bed of
sandstone was also run.  This problem was intended to exercise as many
capabilities of the simulation code as possible, and included a
fluid-poroelastic interface, an orthotropic poroelastic medium with
continuously varying principal axes, and a non-rectilinear
mapped grid designed to conform to the uneven surface of the sandstone
bed.  Results for this demonstration problem were quite complex, with
waves of all three types visible, but the simulation code handled it
without difficulty.

There are many opportunities for extension of the work presented here.
The most obvious route for improvement would be the replacement of the
dimensional splitting scheme with a more accurate method.  While
extending the transverse propagation scheme from~\refmappaper{} into
three dimensions would require an inordinate number of transverse
Riemann solutions, a more promising approach would be to use the
\textsc{SharpClaw} package of Ketcheson et
al.~\cite{ketcheson-parsani-rjl:sharpclaw}, which employs a
semidiscrete approach.  Switching to a semidiscrete scheme would also
allow the use of an exponential
integrator~\cite{hochbruck-ostermann:exponential}, which may allow
better accuracy in the stiff regime identified in~\refcartpaper{}.
Another opportunity to build upon this work would be extension to
higher frequencies --- the numerical scheme used here could be
extended in a straightforward fashion to include additional memory
variables to model a frequency-dependent kernel used to generalize
Darcy's law to higher frequencies~\cite{lu-hanyga:poro-memory-drag}.

In order to facilitate reproduction of these results, all the code
used to produce them has been archived at
\url{http://dx.doi.org/10.6084/m9.figshare.783056}.

\section{Acknowledgements}

This work has benefited greatly from the direct input and advice of
Prof.\ Randall J. LeVeque of the Department of Applied Mathematics,
University of Washington, as well as from the \clawpack{} simulation
framework.
The author also wishes to thank Prof.\ M.\ Yvonne Ou of the
University of Delaware, who introduced him to poroelasticity theory
and whose help has been invaluable in understanding the mechanics of
porous media.
In addition, the treatment of mapped grids in three
dimensions in Section \ref{sec:mapped-grids-3d}, particularly the
appropriate handling of face areas and normals, was inspired by
correspondence with Prof.\ Donna Calhoun of Boise State University.

This work was funded in part by NIH grant 5R01AR53652-2, and by NSF
grants DMS-0914942 and DMS-1216732.

\bibliographystyle{siam}
\bibliography{poro}

\end{document}